\documentclass[asi,article,submit,moreauthors]{Definitions/mdpi} 

\usepackage{caption}
\usepackage{subcaption}
\usepackage{amsmath}
\usepackage{siunitx}
\usepackage{bbold}
\usepackage{vcell}
\usepackage{tabularx}
\usepackage{tabularray}
\usepackage[compatibility=false]{caption}
\usepackage[frozencache=true]{minted}

\usepackage{pgfplots}
\usepgfplotslibrary{colorbrewer}
\usepackage{pgfplotstable}
\usepgfplotslibrary{statistics}
\pgfplotsset{compat = 1.15, cycle list/Set1-8} 
\usetikzlibrary{pgfplots.statistics, pgfplots.colorbrewer} 
\usepackage{pgfplotstable}

\usetikzlibrary{shapes.geometric, arrows, backgrounds, positioning, fit}
\newcommand\ppbb{path picture bounding box}
\tikzstyle{every node}=[font=\small]

\tikzstyle{startstop} = [rectangle, rounded corners, minimum width=3cm, minimum height=1cm, text centered, draw=black, fill=red!30]
\tikzstyle{io} = [trapezium, trapezium left angle=70, trapezium right angle=110, minimum height=1cm, text centered, text width = 2cm, draw=black, fill=blue!30]
\tikzstyle{res} = [trapezium, trapezium left angle=70, trapezium right angle=110, minimum width=3cm, minimum height=1cm, text centered, text width = 2cm, draw=black, fill=green!30]
\tikzstyle{process} = [rectangle, minimum width=3cm, minimum height=1cm, text centered, text width = 2cm, draw=black, fill=orange!30]
\tikzstyle{element} = [rectangle, minimum width=3cm, minimum height=1cm, text centered, text width = 2.5cm, draw=black, fill=black!30]
\tikzstyle{description} = [rectangle, minimum width=3cm, minimum height=1cm, text centered, text width = 3cm, fill opacity=0, text opacity=1, align=left]
\tikzstyle{decision} = [diamond, aspect=1.5, minimum width=3cm, minimum height=1cm, text centered, text width = 2cm, draw=black, fill=white]
\tikzstyle{arrow} = [thick,->,>=stealth]
\tikzstyle{subprocess} = [rectangle, draw=black, fill=orange!30,
                     minimum width=3cm, minimum height=1cm,inner xsep=3mm,
                     text width =\pgfkeysvalueof{/pgf/minimum width}-2*\pgfkeysvalueof{/pgf/inner xsep},
                     align=flush center,
                     path picture={\draw 
    ([xshift =2mm] \ppbb.north west) -- ([xshift= 2mm] \ppbb.south west)
    ([xshift=-2mm] \ppbb.north east) -- ([xshift=-2mm] \ppbb.south east);
                                  }]
\firstpage{1} 
\pubvolume{1}
\issuenum{1}
\articlenumber{0}
\pubyear{2024}
\copyrightyear{2024}
\datereceived{ } 
\daterevised{ } % Comment out if no revised date
\dateaccepted{ } 
\datepublished{ } 
\hreflink{https://doi.org/} % If needed use \linebreak
\Title{An Application Driven Method for Assembling Numerical Schemes for the Solution of Complex Multiphysics Problems}

\TitleCitation{An Application Driven Method for Assembling Numerical Schemes for the Solution of Complex Multiphysics Problems}

\Author{Patrick Zimbrod $^{1,\dagger}$\orcidA{}, Michael Fleck $^{2}$\orcidB{} and Johannes Schilp $^{1}$*}

\AuthorNames{Patrick Zimbrod, Michael Fleck and Johannes Schilp}

\AuthorCitation{Zimbrod, P.; Fleck, M.; Schilp, J.}
\address{%
$^{1}$ \quad Applied Computer Science, University of Augsburg\\
$^{2}$ \quad Metals and Alloys, University of Bayreuth}

\corres{Correspondence: patrick.zimbrod@uni-a.de}

\abstract{Within recent years, considerable progress has been made regarding high-performance solvers for Partial Differential Equations (PDEs), yielding potential gains in efficiency compared to industry standard tools. However, the latter largely remains the status quo for scientists and engineers focusing on applying simulation tools to specific problems in practice.
We attribute this growing technical gap to the increasing complexity and knowledge required to pick and assemble state-of-the-art methods.
Thus, with this work, we initiate an effort to build a common taxonomy for the most popular grid-based approximation schemes to draw comparisons regarding accuracy and computational efficiency.
We then build upon this foundation and introduce a method to systematically guide an application expert through classifying a given PDE problem setting and identifying a suitable numerical scheme.
Great care is taken to ensure that making a choice this way is unambiguous, i.e. the goal is to obtain a clear and reproducible recommendation.
Our method not only helps to identify and assemble suitable schemes but enables the unique combination of multiple methods on a per-field basis.
We demonstrate this process and its effectiveness using different model problems, each comparing the resulting numerical scheme from our method with the next best choice.
For both the Allen Cahn and advection equations, we show that substantial computational gains can be attained for the recommended numerical methods regarding accuracy and efficiency.
Lastly, we outline how one can systematically analyze and classify a coupled multiphysics problem of considerable complexity with 6 different unknown quantities, yielding an efficient, mixed discretization that in configuration compares well to high-performance implementations from the literature.}

\keyword{simulation; multiphysics; finite difference; finite volume; finite element; discontinuous galerkin}

\begin{document}
% For arXiv
\nolinenumbers

%%%%%%%%%%%%%%%%%%%%%%%%%%%%%%%%%%%%%%%%%%
%\endnote{This is an endnote.} % To use endnotes, please un-comment \printendnotes below (before References). Only journal Laws uses \footnote.

% The order of the section titles is different for some journals. Please refer to the "Instructions for Authors” on the journal homepage.

\section{Introduction}\label{sec:introduction}

Within the discipline of accurately and efficiently modeling physical processes that are described by Partial Differential Equations (PDEs), one is confronted with an increasing amount of choices of numerical methods to perform this job.
To an application-oriented expert, that is, engineers as well as scientists who are quite familiar with their respective scientific domain, but not necessarily with the numerics of PDE approximation, this abundance of choices may quickly appear daunting due to the growing amount of research in the numerics community.
This trend can be observed in particular within the Finite Element framework, which has grown to be one of the most widely used and also formally investigated methods in this field.
Modern, general and mathematically rigorous implementations of the Finite Element Method (FEM), such as \texttt{deal.II}\,\cite{bangerthDealIIGeneralpurpose2007}, \texttt{FEniCS}\,\cite{AlnaesEtal2015}, \texttt{Firedrake}\,\cite{FiredrakeUserManual} or \texttt{MFEM}\,\cite{andersonMFEMModularFinite2021a}, offer a wide variety of formulations that one can choose from, with varying degrees of customizability.
The amount of different options is best illustrated by considering the various Finite Elements and corresponding function spaces that the application expert typically can and must choose from nowadays.
For instance, the popular website DefElement, which summarises a vast amount of Finite Element types along their characteristics and shape functions, offers over 45 different choices of element to approximate scalar and over 40 for vector quantities\,\cite{scroggsDefElementEncyclopediaFinite2023}.
These elements are, for the most part, readily implemented in the abovementioned software libraries.
However, not all of them necessarily produce good or even stable approximations for any given PDE problem \cite{brezziDiscourseStabilityConditions1990}.
Even if one were to have prior knowledge on a good choice of function space, e.g. $H(\textrm{div})$ to approximate divergence-free quantities\,\cite{johnDivergenceConstraintMixed2017}, one would still have to choose between more than 20 different kinds of Finite element.
This exemplifies the large flexibility and sophistication of the FEM that has evolved in recent years, but consequently also the ambiguity of choice that these developments bring along.
The outlined problem of choosing the right tool for the problem naturally becomes even more of an issue when other standard tools that are used in practice are additionally considered, such as Finite Difference and Finite Volume methods.
These recent advancements are, however, in contrast to established numerical techniques that are oftentimes still used as a standard in practice.
For example, the Finite Volume method, although relying on a rather old and simple concept, is still widely considered the standard practice for solving problems in Computational Fluid Dynamics\,\cite{shademanEvaluationOpenFOAMAcademic2013}.
In contrast, there have been several recent developments based on the Discontinuous Galerkin method that have shown to perform noticeably better for fluid dynamics problems \cite{zhouNumericalComparisonWENO2001}.
We thus note that a gap has emerged between industry standard tools and modern, high-performance numerical methods that we attribute to the increasing complexity and insight required to properly assemble the latter.

Additionally, apart from these outlined recent developments, choosing an approximation method that offers the right amount and type of degrees of freedom is essential for obtaining a solution that can be efficiently computed.
That is, a numerical scheme should be tailored to a specific problem at hand to achieve good convergence and stability properties.
Where the former might be negligible in practice since achieving the utmost performance is not always important, having a stable approximation is paramount.
One must hence make use of the specific properties of a PDE problem to yield a good approximation.
This line of reasoning not only applies to problems governed by a single PDE but even more so to multiple equations forming a system.
We will in the following denote such problems as multiphysics problems.
These oftentimes describe processes that are either distinct in their qualitative behavior or even operate on different length scales.
As a result, using one discretization method in a monolithic way will not perform equally well for each PDE of that system, creating bottlenecks regarding either stability or accuracy and thus further complicating the question of which specific method is best suited for the job.
On the other hand, using multiple numerical methods within one multiphysics problem requires interoperability of all discretizations.
Coupling different solvers has shown to quickly become tedious regarding implementation and may even yield severe bottlenecks due to large amounts of data transfer.
Having a common formulation for some schemes is thus desirable, such that one may recover specific methods by imposing abstractions, for example by omitting some steps in assembling a global linear system.

In this work, we make an initial effort towards closing the outlined, increasing gap between state-of-the-art research in the numerics community and best practices in applications.
As covering the entirety of numerical methods available is impossible in practice, especially for non-experts in every single aspect, we initially restrict the scope of view to a rather narrow subset of methods.
The added benefit of this approach is that this enables us to remove methods that would otherwise render making a problem-oriented choice largely ambiguous.

The contributions of this work to address the abovementioned problems are twofold:
First, we propose a unifying approximation taxonomy that enables recovering the most prevalent grid-based numerical methods by imposing some well-defined abstractions, albeit with a relatively narrow scope for now.
Secondly, we propose a generalized framework for choosing an appropriate, that is, stable and performant numerical scheme given a fixed set of inputs.
These include the system of PDEs, the triangulation where the problem is defined as well as the available computing hardware.
As such a method necessitates a unified view of all the numerical methods considered to enable quantifiable comparisons, this heavily builds on the first part of this article.
In combination, this enables an application expert to make an informed choice on which scheme to use on a per-field basis given some equally well-defined inputs.
Since in an application setting, stability is typically more important than convergence rate, we focus our effort on providing methods that produce stable solutions but do not fall behind too much compared to the utmost performant alternatives.

The remainder of this work is thus structured as follows:
We first present a brief review of works in the literature that are concerned with comparing the mentioned numerical schemes which we subsequently build upon.
The theory necessary to construct such a baseline will be covered in the following sections.
We will then proceed by analyzing typically given inputs for a PDE problem and assemble a method to systematically derive suitable numerical schemes.
The results of choosing and implementing numerical methods according to our developed framework will be demonstrated afterward.
We will investigate two distinct benchmark problems, each comparing two methods that would be closest to being optimal in that specific case.
Special attention will be given to accuracy concerning the analytical solution as well as computational complexity and performance.
Finally, we will demonstrate the effectiveness of the presented framework by assessing a complex and currently relevant multiphysics problem.
We show how to systematically arrive at a mixed choice of discretization schemes using the proposed decision method. Finally, we indicate some relevant works in the literature that employ similar approximations, highlighting the validity and relevance of this work.

\section{Previous Works}

In this section, we outline some prior efforts aimed at comparing different grid-based numerical methods or drawing connections between them.
Besides the works mentioned henceforth, there exists a vast amount of literature that is concerned with the unique properties of each method, which we omit here for the sake of brevity.

In the literature, there are rather few works that are concerned with spanning the connection between different grid-based approximation schemes.
Some authors rigorously showed the equivalence of the Finite Volume method to either Mixed Finite Element \cite{barangerConnectionFiniteVolume1996} or Petrov Galerkin Finite Element Methods \cite{yeRelationshipFiniteVolume2001,idelsohnFiniteVolumesFinite1994}.
With regards to the Finite Difference and Finite Element Method, Thomèe has shown early on that the FEM can be understood as a somewhat equivalent, yet generalized variant of taking Finite Differences on arbitrary grids \cite{thomeeFiniteDifferenceFinite1984}.
Some general differences between these schemes were outlined by Key and Krieg \cite{keyComparisonFiniteElementFiniteDifference1973}.
In the work of Shu, some analogies have been brought up between Finite Volume and Finite Difference schemes in WENO formulation \cite{shuHighorderFiniteDifference2003}.
A theoretical and numerical comparison between higher-order Finite Volume and Discontinuous Galerkin Methods was conducted by Zhou et al. \cite{zhouNumericalComparisonWENO2001}.
Additionally, Dumbser et al. constructed a unifying framework to accommodate high-order Finite Volume and Discontinuous Galerkin schemes \cite{dumbserUnifiedFrameworkConstruction2008a}.
In the context of elliptic PDEs, Lin et al. present a theoretical and empirical comparison between the comparably new weak Galerkin, Discontinuous Galerkin and mixed Finite Element schemes \cite{linComparativeStudyWeak2015a}.
A comparative study between Discontinuous Galerkin and the Streamline Upwind Petrov Galerkin method for flow problems can be found in \cite{yurunComparativeStudyDiscontinuous1997}.
These works in summary draw point-wise comparisons between some grid-based approximation schemes.
Despite being quite useful for disseminating individual advantages and disadvantages for a given application, one may still lack an understanding of the general properties.
Furthermore, Bui-Thanh has presented an encompassing analysis and application of the Hybridizable Discontinuous Galerkin Method (HDG) to solve a wide variety of PDE-governed problems.
It was therefore shown that this numerical scheme is general and powerful enough to form a unified baseline \cite{bui-thanhGodunovUnifiedHybridized2015}.
In addition, due to the generality of this method, there have been works that attempt to benchmark DG methods to more conventional and widely adopted Continuous Galerkin methods (CG) \cite{yakovlevCGHDGComparative2016,kronbichlerPerformanceComparisonContinuous2018,kirbyCGHDGComparative2012}.
Some authors proposed combinations of numerical schemes that operate optimally to solve hyperbolic \cite{gaburroUnifiedFrameworkSolution2021} or parabolic \cite{yangOptimallyEfficientTechnique2017} systems of PDEs.

In summary, we draw the following conclusion from this brief review of the relevant literature.
To this date, there only exist a few comparisons between grid-based approximation schemes that outline common properties in a mostly ad hoc or point-wise manner.
We have, however, outlined in the previous section that it would be beneficial from an application-oriented perspective to have an encompassing taxonomy for these numerical schemes.
Furthermore, many different variants of numerical schemes have been proposed to tackle a wide variety of PDE problems.
What appears to be missing, though, is a general guideline on how to choose between these vast alternatives to obtain a method, or possibly a combination of different methods, to solve a system in a stable (above all) and reasonably efficient way.

\section{Theoretical Baseline}\label{sec:baseline}

The general procedure of this chapter is as follows.
We first introduce the most general scheme considered here, which is the Discontinuous Galerkin Method.
Then, for each additional scheme considered, we individually work out the necessary simplifications to arrive at that numerical method starting at the DGM.
In the remaining sections of this article, we use underline notation ($\underline{\cdot},\,\underline{\underline{\cdot}}$) to indicate vectors and matrices and roman indices ($i,j$) to denote elements of lists or arrays on the computational level.
We assume the reader of this article to have a coarse overview of the presented methods, but not much insight into the specifics of each. We thus present the necessary theory in a comparably coarse manner that focuses on the qualitative characteristics.

\subsection{Discontinuous Galerkin Method}

This method was originally proposed in 1973 to solve challenging hyperbolic transport equations in nuclear physics \cite{reedTriangularMeshMethods1973}.
In spirit, it can be held as a synthesis of Finite Element and Finite Volume schemes and poses a generalized variant of both.

To derive such a scheme, we begin by stating the strong form of a given PDE. The most straightforward example in this case would be a first-order linear advection equation with a homogeneous von Neumann boundary condition:
\begin{equation} \label{eq:advection-strong}
    \partial_t \alpha + \underline{u} \cdot \nabla \alpha = 0, \qquad \frac{\partial \alpha}{\partial n} = 0 \quad\forall x \in \partial \Omega
\end{equation}
where $\partial_t$ signifies the temporal derivative, $x$ is the set of spatial coordinates, $\partial \Omega$ denotes the boundary of the computational domain $\Omega$ and $n$ is the unit normal with respect to $\partial \Omega$.
In this case, and henceforth in this article, we assume for reasons of simplicity that the velocity field is divergence-free, i.e. $\nabla \cdot \underline{u} = 0$.

We now state the weak form of Eq. \ref{eq:advection-strong}, that is we multiply with a test function $v$, integrate over the entirety of the domain $\Omega$, and apply partial integration to the second term on the left-hand side that contains the nabla operator.
For a divergence free velocity field, one may set $\underline{u} \cdot \nabla \alpha = \nabla \cdot \left(\underline{u}\alpha\right)$, which results in the following formulation: Find $\alpha \in V$ such that:
\begin{equation}\label{eq:advection-weak-dgm}
    \int_\Omega v \cdot \partial_t \alpha \,dx + \int_{\partial \Omega} v \alpha (\underline{u}\cdot \underline{n}) \,ds - \int_\Omega \nabla v \cdot \underline{u} \alpha \,dx = 0 \qquad \forall v \in W
\end{equation}
Where we need to make an appropriate choice for the solution space $V$ and the test space $W$ which may but do not need to differ from each other.

Due to partial integration, we now encounter an additional term that has to be integrated over the domain boundary $\partial \Omega$, where $(\underline{u}\cdot \underline{n})$ denotes the velocity component normal to the boundary.

To make such a problem solvable by a computer, one must additionally choose the discretization of the solution space $V$, denoted $V_h$.
A particularly popular choice of space is the set of Lagrange polynomials.
In addition, the physical space must be discretized in the form of a triangulation.
The DG scheme then consists of assembling the finite-dimensional, linear system on the element level.
This enables high locality of the solution process, which leads to efficient computation on parallel architectures as less data transfer is required.

One resulting key feature of the DG scheme is that the elements now do not overlap anymore in terms of their degrees of freedom.
Thus, the global problem is broken up into individual problems.
This in general leads to large systems that are however sparse and in the case of the mass matrix even block diagonal.
The remaining term, often denoted the numerical flux, is the only term within the physical domain that ensures coupling across elements.
Through evaluation of this surface integral, adjacent degrees of freedom are coupled and thus global conservation of quantities can be assured.

As the polynomial space of DG schemes only belongs to the $L^2$ space of functions but not $H^1$, the basis functions are discontinuous and thus the derivative at boundaries is not well defined. 
Solving PDEs involving second derivatives is thus not possible as is.
As a consequence, there have been many successful extensions of this method to circumvent that problem. At this point, we name the most prevalent schemes, namely the Symmetric Interior Penalty \cite{epshteynSymmetricInteriorPenalty}, Hybridizable \cite{warburtonDiscontinuousGalerkinMethod1999} and Local DG scheme \cite{cockburnLocalDiscontinuousGalerkin1998}.
These methods, despite having different approaches, have been extensively studied and compared to each other \cite{arnoldUnifiedAnalysisDiscontinuous2002}.
As it turns out, all methods work well and have individual advantages and disadvantages.
For this work, we will take the Hybridizable DG scheme as a general framework.
We note here that the proposed method would however work with any of the other schemes given above.

Within the abovementioned methods, one introduces an additional term in the weak form that serves as a penalty for discontinuous solutions.
An alternative approach that is also pursued within the HDG scheme is the algebraic manipulation of the PDE system by splitting.
One recursively introduces new dependent variables for quantities that appear in higher-order derivatives such that each quantity is differentiated at most once.
We illustrate this using the Poisson equation:
\begin{equation}
    \Delta u = 0
\end{equation}
The corresponding, well-known weak form is: Find $u \in V$, such that for all $v \in V$
\begin{equation}\label{eq:poisson-weakform}
    \int_{\partial\Omega} v (\nabla u \cdot \underline{n}) \,ds - \int_\Omega \nabla v \nabla u \,dx = 0
\end{equation}
By introducing the auxiliary variable $\sigma = \nabla u$, Eq. \ref{eq:poisson-weakform} is extended to the following system:
\begin{align}
    \int_{\partial\Omega} v (\sigma \cdot \underline{n}) \,ds - \int_\Omega \nabla v \sigma \,dx = 0 \\
    \sigma = \nabla u
\end{align}
Employing such an approach enables splitting PDE systems of arbitrary order resulting in larger systems of first order PDEs.

\subsection{Continuous Galerkin Finite Element Method}

The most straightforward step to conduct is to derive the Continuous Galerkin (CG) from the DG method.
The former is oftentimes also referred to as the classic Finite Element method, being the original formulation used to solve problems in structural mechanics \cite{liuEightyYearsFinite2022}.

In this case, all degrees of freedom (DoFs) in the domain are global, in contrast to being local to each cell.
However, each basis function associated with a given degree of freedom has compact support and is thus only non-zero within the direct vicinity.
The resulting linear system hence remains sparse but has considerably fewer DoFs than an equivalent discretization produced by a DG method.

One may obtain a CG method starting from the DGM by strongly coupling the degrees of freedom at cell interfaces. In other words, the previously discontinuous approximation must be made continuous.
In terms of the weak form of a given problem, the numerical flux that has been introduced by partial integration has to vanish.
This step is exactly taken in deriving weak forms for the CG scheme.
The equivalent weak form of the advection equation given by Eq. \ref{eq:advection-weak-dgm} is then: Find $u \in V$ such that
\begin{equation}\label{eq:advection-weak-cgm}
    \int_\Omega v \cdot \partial_t \alpha \,dx - \int_\Omega \nabla v \cdot \underline{u} \alpha \,dx = 0 \qquad \forall v \in V
\end{equation}
By coupling coinciding DoFs, one may equivalently introduce shared DoFs between cells.
This results in comparison to the DGM in a smaller global system that is in turn more coupled, yielding more non-zero entries per row and column in the system matrices.

\begin{figure}[htbp!]
    \centering
    \captionsetup[subfigure]{justification=centering}
    \begin{subfigure}[c]{0.49\textwidth}
        \centering
        \includegraphics[]{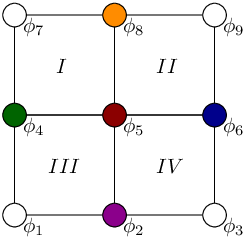}
        \caption{}
        \vfill
        \label{fig:cg-2d}
    \end{subfigure}
    \begin{subfigure}[c]{0.49\textwidth}
        \centering
        \includegraphics[]{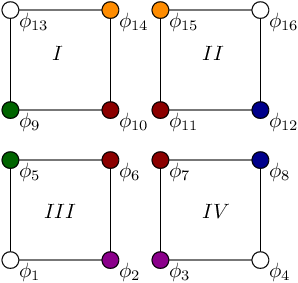}
        \caption{}
        \label{fig:dg-2d}
    \end{subfigure}
   \caption{Coupling of global DoFs in the Continuous Galerkin (a) versus Discontinuous Galerkin FEM (b), both of first order. In the latter case, DoFs are entirely local to the cell and thus receive no contribution from neighboring cells. Weak coupling is only introduced by the additional numerical flux. Coupled DoFs are drawn in identical colors.}
   \label{fig:cg-vs-dg-2d}
\end{figure}

The condensation of such a system by coupling DoFs is illustrated in Figure \ref{fig:cg-vs-dg-2d}.
From the numbering of DoFs in both figures, it becomes apparent that the amount of additional allocations grows drastically with increasing dimensionality of the problem.

As the numerical flux is zero by definition for a CG scheme, we may also omit it from computation.
Thus, the CGM is noticeably less arithmetically intensive in this regard.
However, this computational saving is offset by the strong coupling of DoFs, resulting in a more dense linear system and possibly a more complex assembly process in terms of memory management.

As equivalence can be shown here based on the weak form and thus early on in the model assembly process, the choice of Finite Element is unaffected.
This in consequence also applies to the chosen type of triangulation or the order of approximation.

\subsection{Finite Difference Method}\label{sec:fdm}

\begin{figure}[htbp!]
    \centering
    \captionsetup[subfigure]{justification=centering}
    \begin{subfigure}[t]{0.49\textwidth}
        \centering
        \includegraphics[]{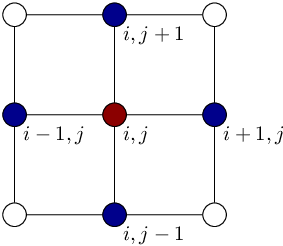}
        \caption{}
        \label{fig:fd-stencil-2d}
    \end{subfigure}
    \hfill
    \begin{subfigure}[t]{0.49\textwidth}
        \centering
        \includegraphics[]{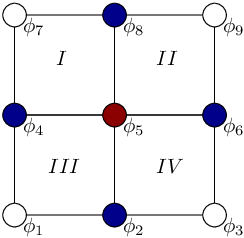}
        \caption{}
        \label{fig:cg-2d-collocatoin}
    \end{subfigure}
   \caption{Comparison of the nodal nature of the FDM (a) versus the cell-wise assembly used in the CG FEM (b) for an identical, cartesian triangulation with 9 nodes. Both methods are formulated as first-order approximations.
   Red-colored nodes signify the points where the PDE is evaluated. Contributions to this node are taken from blue nodes, whereas white nodes from no contribution.}
   \label{fig:fd-vs-cg-2d}
\end{figure}

At first glance, the Finite Difference Method (FDM) might appear to be conceptually different from the Finite Element methods given above. Instead of treating the discretized problem in an element-wise manner, the FDM operates on discrete points directly and per se lacks a notion of cells in the domain.
Yet, both methods still may yield identical results in discretization.
A comparison of both approaches is shown in Figure \ref{fig:fd-vs-cg-2d}.

In the figures, $i$ and $j$ denote vertical and horizontal indices of grid nodes, $\phi_i$ are the FE basis functions and Roman letters denote indices of cells.
Forming for example the laplacian for an FDM requires access to the vicinity of vertex $i,j$ (red node) in all cartesian directions (blue-colored nodes). For both methods, gray-marked nodes do not pose a contribution to the value of the central node. For a special case of FEM with quadrilateral elements, the same nodes form contributions to the global basis function $\phi_5$. However, we now do not evaluate the laplacian operator directly but instead gather contributions from weak form integrals. In the case of $\phi_5$, we have to gather contributions from cells $I$ to $IV$.

However, one may still show the equivalence of CGM and FDM by investigating the resulting global linear system.
We exemplify this claim using the laplacian as a differential operator and the second-order central stencil.
This formulation continues to be widely used as an approximation technique.
As will be shown later, this particular choice of operator is especially straightforward to compare with the CG-FEM due to the choice of trial and test functions.
On a cartesian, two-dimensional grid with uniform spacing $h$ in both directions, the approximation reads
\begin{equation}\label{eq:fd-laplacian}
    \Delta u \approx \frac{u_{i-1,j}+u_{i,j-1}-4u_{i,j}+u_{i+1,j}+u_{i,j+1}}{h^2}.
\end{equation}
Such a system in stencil notation will produce a global matrix with main diagonal values 4 and four off-diagonals with entries 1.

We now proceed to construct an equivalent CG Finite Element scheme, where the global system matrix is required to be exactly equivalent to the FD formulation.

The weak (CG) laplacian can be formulated as: Find $u_h \in V$ such that for all $v_h \in V$
\begin{equation}
    \int_\Omega v \Delta u \,dx = \underbrace{\int_{\partial\Omega} v (\nabla u \cdot \underline{n}) \,ds}_{=0} - \int_\Omega \nabla v \nabla u \,dx
\end{equation}
We have in this case introduced the additional restriction that trial and test space be identical, that is, we use a Bubnov Galerkin method.
Now, let $\Omega$ be an identical triangulation to the FD variant using quadrilateral $\mathbb{Q}^1$ elements, that is linear Lagrange elements.

Then, the four basis functions spanning the reference element are:
\begin{align}
    \phi_1(x,y) &= xy - x - y + 1 &&\\
    \phi_2(x,y) &= x (1 - y) &&\\
    \phi_3(x,y) &= y (1 - x) &&\\
    \phi_4(x,y) &= x y &&
\end{align}
The Finite Difference stencil given by Eq. \ref{eq:fd-laplacian} only takes into account contributions from nodes that lie strictly horizontally or vertically from the node of interest.
As a consequence, the node on the reference quadrilateral that is positioned diagonally from the center node must not have any contribution to the weak form integral, otherwise the resulting linear system cannot be equal.
We thus need to evaluate the weak form in a way such that the resulting matrix $\phi_i \cdot \phi_j$ becomes sparse.
It turns out that this can be achieved by choosing a collocation method for quadrature.
In that case, quadrature points are chosen to coincide with the node coordinates and as a consequence, the mass matrix $\phi_i \cdot \phi_j$ becomes the identity matrix.

From the family of Gaussian quadrature schemes, one can achieve this using a Gauss-Lobatto quadrature of order equal to the polynomial order of the Finite Element.
We note at this point that choosing this particular combination of quadrature method and number of nodes will lead to inexact integration and thus, a numerical error is introduced.
In this case, to produce a collocated scheme, one must pick the second-order Gauss-Lobatto variant using two quadrature nodes per coordinate direction.
As this type of integration is known to be accurate up to degree $2n-3$, this scheme will only integrate linear polynomials exactly.
However, the above-listed basis functions are bilinear and have a combined polynomial order of $2$.
Thus, integration will not be accurate in this case.
However, the modern FEM in general does not prescribe any particular method of integrating the weak formulation per se \cite{ciarletFiniteElementMethod2002}. 
Thus, although the results of a properly implemented FEM in the sense of exact integration will be slightly more accurate, one can still show equivalence regarding a particular instance of the FEM.

We now evaluate the element-wise stiffness matrix $-\int_{\Omega^{(e)}} \nabla_k \phi_i \nabla_k \phi_j \,dx$ within the reference domain $\left[0;1\right]\times\left[0;1\right]$ for the given first order Lagrange element using Gauss-Lobatto quadrature, more specifically the variant using two quadrature points per coordinate direction.
This results in:
\begin{equation}
    \underline{\underline{K}}^{(e)} = \begin{bmatrix}
         -1  & 1/2  & 1/2  & 0 \\
         1/2 &  -1  & 0    &  1/2 \\
         1/2 & 0    & -1   &  1/2 \\
        0    &  1/2 &  1/2 & -1
    \end{bmatrix}
\end{equation}
As such, $\underline{\underline{K}}^{(e)}$ does not yet equal Eq. \ref{eq:fd-laplacian}.
The final step consists of assembling the linear system in the physical domain using the reference stiffness matrix.
In a cartesian mesh in two dimensions, an interior node is owned by four quadrilateral elements.
If one carries out this assembly process an equivalent formulation can be obtained:
\begin{equation}
    \underline{\underline{K}} = \begin{bmatrix}
        \ddots & \ddots &  &  & \vdots & \vdots & \vdots &  &  &  &  \\
        \ddots & \ddots & \ddots &  & \vdots & 1 & \vdots &  &  &  &  \\
         & \ddots & \ddots & \ddots & \vdots & \vdots & \vdots &  &  &  &  \\
         &  & \ddots & \ddots & \ddots & 1 & \vdots &  &  &  &  \\
        \cdots & \cdots & \cdots & \ddots & \ddots & \ddots & \vdots & \cdots & \cdots & \cdots & \cdots \\
        \cdots & 1 & \cdots & 1 & \ddots & -4 & \ddots & 1 & \cdots & 1 & \cdots \\
        \cdots & \cdots & \cdots & \cdots & \vdots & \ddots & \ddots & \ddots & \cdots & \cdots & \cdots \\
         &  &  &  & \vdots & 1 & \ddots & \ddots & \ddots &  &  \\
         &  &  &  & \vdots & \vdots & \vdots & \ddots & \ddots & \ddots &  \\
         &  &  &  & \vdots & 1 & \vdots &  & \ddots & \ddots & \ddots \\
         &  &  &  & \vdots & \vdots & \vdots &  &  & \ddots & \ddots 
        \end{bmatrix}
\end{equation}
The exact position of the one entry in the typically large and sparse matrix depends on the mesh topology as well as the global numbering of the degrees of freedom.

For the discretization of other operators, a similar argument holds, as the shown procedure is irrespective of the choice of weak form or basis function.
For example, one could discretize the gradient of a function $\nabla u$ using an upwind Finite Difference formulation in fluid mechanics for resolving convective terms.
An equivalent Finite Element method can be assembled by producing a weak form as given in the above example, choosing the same collocation method and carrying out the integration numerically.
However, one important difference is that one cannot choose the test space to be equivalent to the trial space.
This would yield a symmetric system that does not correspond to an upwind Finite Difference formulation and is also not stable in the case of solving a pure advection equation.
One must instead choose a test space with asymmetric test functions to account for the notion of an upwind node, thus yielding a Petrov Galerkin scheme \cite{brooksStreamlineUpwindPetrovGalerkin1982}.

We can as a result summarise the FDM to be a special instance of the CG FEM.
On the one hand, integration is restricted to a collocation method and on the other hand the Jacobian mapping from reference to physical elements is constant throughout the domain.
This close relationship has also been hinted at by analysis of boundary value problems by Thomèe \cite{thomeeFiniteDifferenceFinite1984}.

For the sake of achieving the same discretization, the use of Finite Differences over Finite Elements becomes apparent from the discussion above.
Most strikingly, the process of producing a local stencil is vastly more straightforward than performing element-wise assembly and gathering the weak form integrals in a global, sparse linear system.
Each element-wise operation in assembly would otherwise require the evaluation of the mesh jacobian for the requested element, that is the mapping from the reference to physical space.
Furthermore, this constant stencil enables Finite Difference schemes to operate in a matrix-free manner easily.
For larger systems, this can help to avoid a large amount of allocated memory, thus being suited well for modern hardware architectures that are typically memory-bound.

These advantages are however offset by some topological restrictions on the mesh.
The simplicity of a constant stencil also implies that the mesh must not deviate from a cartesian geometry.
Otherwise, additional complexity is introduced since Eq. \ref{eq:fd-laplacian} becomes a stencil in the reference domain that has to be mapped to the physical domain.
This would still save the computational effort to assemble the weak form.
However, since this process only has to be carried out once for the reference element, the computational impact can be held low by pre-computing the integrand.

\subsection{Finite Volume Method}\label{sec:fvm}

In a similar vein to the FDM, the use of Finite Volumes might appear distinctly different from the idea of Finite Elements.
Here, we make extensive use of Stokes' theorem to replace volume with hull integrals in conservation laws \cite{eymardFiniteVolumeMethods2000}.
There exist different formulations of this method, most namely a cell- and vertex-centered form.
The main difference lies in where the solution is stored.
In the former case, the solution is stored at the polygonal cell centers that are spanned by the mesh vertices.
The latter, instead, directly uses these vertices as solution points \cite{moukalledFiniteVolumeMethod2016}.
In this case, one does not operate on the computational mesh directly but rather on its dual.
As the cell-centered formulation is considerably more widely used, we investigate this variant further in the following.

It can be shown however that the FVM can simply be considered a Bubnov Discontinuous Galerkin method of polynomial order zero.
To illustrate this, we again turn to Eqs. \ref{eq:advection-strong} and \ref{eq:advection-weak-dgm} describing the strong and weak form of the advection equation.
A Finite Volume approximation in conservation form is:
\begin{equation}\label{eq:advection-fvm}
    \int_\Omega \partial_t \alpha \,dx + \int_\Omega \underline{u} \nabla \alpha \,dx =  \int_\Omega \partial_t \alpha \,dx + \int_{\partial\Omega} \alpha (\underline{u}\cdot\underline{n}) \,ds
\end{equation}
Apart from the presence of a test function $v$ in Eq. \ref{eq:advection-weak-dgm}, the second integrand simply represents the net flux of the conserved quantity $\underline{u}\alpha$ over the set of element boundaries.

For Eqs. \ref{eq:advection-weak-dgm} and \ref{eq:advection-fvm} to be equivalent in this case, the third integrand resulting from partial integration has to vanish in addition.
However, this can be shown trivially by setting the order of the polynomial space for the trial and test function to zero.
Then, the derivative of the test function vanishes and thus the entire term does not contribute to the weak form.

\begin{figure}[htbp!]
    \centering
    \captionsetup[subfigure]{justification=centering}
    \begin{subfigure}[c]{0.49\textwidth}
        \centering
        \includegraphics[]{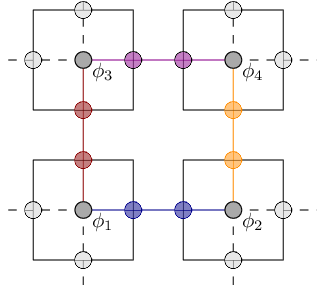}
        \caption{}
        \label{fig:fv-2d}
    \end{subfigure}
    \hfill
    \begin{subfigure}[c]{0.49\textwidth}
        \centering
        \includegraphics[]{Figs/Theory/dg-2d.eps}
        \caption{}
        \label{fig:dg-2d_2}
    \end{subfigure}
   \caption{Comparison of FVM (a) versus DGM (b) on an identical quadrilateral triangulation.
   Coupled DoFs are marked in identical colors.
   For the FVM, one must first reconstruct the values of the DoFs at the mesh facets to then compute the numerical flux.}
   \label{fig:fv-vs-dg-2d}
\end{figure}

After performing this step, the test function is still present in the remaining parts of Eq. \ref{eq:advection-weak-dgm}.
For the remaining terms to be equivalent, they must vanish out of the equation as well.
This can be accomplished straightforwardly by fixing the value of the test function to be unity.
In the weak form, this step is admissible since it must hold for all instances of $V$.
As $1 \in V^0$ where $V^0$ is the space of constant polynomials, this statement holds in particular for a Bubnov Galerkin scheme, as trial and test space must be identical.
The qualitative similarity of both schemes is illustrated in Figure \ref{fig:fv-vs-dg-2d}.

For both schemes, DoFs are entirely local to the cell and coupling happens through the calculation of a numerical flux - or in more formal terms, through the evaluation of the hull integral in the corresponding weak form.
However, the FVM only stores one DoF per cell which has notable implications for the calculation of the numerical flux.
This means that as a first step, the cell values have to be reconstructed at the mesh facets. These reconstructed DoFs which depend on the cell values that they interpolate between can then be coupled to their counterparts at opposing mesh facets. These relationships are denoted by DoFs being colored identically (coupled) and being transparent (reconstructed) in Figure \ref{fig:fv-2d}.
For a DGM scheme of first order or higher, this interpolation step oftentimes is not necessary if the quadrature scheme is chosen carefully.
For a collocation method (see section \ref{sec:fdm} for a more thorough discussion), one does not need to tabulate the full list of DoF values at the set of facet quadrature points, but rather only a small subset of DoFs that are owned by the facet \cite{hesthavenNodalDiscontinuousGalerkin2008}.

In summary, the FVM can again be considered as a special instance of the Bubnov type DGM, where the shared polynomial space $V$ is taken to be of constant order and the test function $v$ is set to be unity.

We note similarly to the discussion on the FDM that this simplification of Eq. \ref{eq:advection-weak-dgm} brings with it some computational advantages that can be offset by sacrificing flexibility.
The absence of a true weak form in a Finite Volume formulation again means that actual assembly is not needed.
In addition, one may omit the transformation from the reference to physical space, as interpolating degrees of freedom to mesh facets and forming a finite sum of these contributions can be done on the mesh directly.
The caveat of this approach is that FVM in principle is bound to be at most first-order accurate.
In practice, this does not hold as the FVM can be extended to higher orders by applying higher order flux reconstruction techniques \cite{zhouNumericalComparisonWENO2001,shuHighorderFiniteDifference2003}.
Such techniques can however quickly become computationally expensive as well with increasing order.
This is achieved in this case by widening the stencil for polynomial reconstruction, increasing memory and time complexity by a considerable amount \cite{liuRobustReconstructionUnstructured2013}.

\subsection{Summary}

In this section, we have established a common framework to formulate the most prevalent grid-based numerical schemes for the solution of PDEs.

It turns out that the DG Method possesses enough flexibility to incorporate the CGM, FDM and FVM by imposing a set of restrictions.
A summary of the results presented in this section on how the schemes compare overall is given in Table \ref{tab:schemes-comparison}.

\begin{table}
    \centering
    \caption{Comparison of the individual restrictions that the presented schemes impose. Certain simplifications bring with them computational advantages, as discussed above.}
    \label{tab:schemes-comparison}
    \resizebox{\linewidth}{!}{%
    \begin{tblr}{
      hline{1,6} = {-}{0.08em},
      hline{2} = {-}{0.05em},
    }
    \textbf{Scheme} & \textbf{Geometry} & \textbf{Function Space} & \textbf{Weak Form} & \textbf{Quadrature}\\
    DGM & Arbitrary & Discontinuous ($L^2$) & Full & Arbitrary\\
    CGM & Arbitrary & Continuous ($H^1$) & {No hull integrals \\over interior facets} & Arbitrary\\
    FDM & Cartesian & Continuous ($H^1$) & {No hull integrals \\over interior facets} & Collocation\\
    FVM & Arbitrary & {Discontinuous ($\mathbb{P}^0 \in L^2$),\\Bubnov Galerkin} & No volume integrals & None required
    \end{tblr}
    }
\end{table}

This framework is not only of theoretical use.
Rather, such a common formulation also enables us to combine these schemes arbitrarily to solve larger problems.
As each scheme possesses strengths and means to gain computational efficiency, this is an important result since it enables efficiently mixed discretizations of multiphysics problems.
Establishing a practical method to achieve exactly this will be the content of the next section.

Before concluding the discussion on relating the above numerical schemes, we add an important remark.
There do exist several extensions to these methods that in general do not fit into the framework that has been established.
We will list a few examples for the sake of illustration.

There do exist formulations of the FDM that can capture domains with less regularity, see for example \cite{fornbergGenerationFiniteDifference1988,visbalUseHigherOrderFiniteDifference2002,zhangThreedimensionalElasticWave2012,perroneGeneralFiniteDifference1975}.
One can also find alternative discretization methods based on FDM in the literature that encompass the notion of missing structure in grids more naturally, such as understanding vertices as centroids of Voronoi cells \cite{sukumarVoronoiCellFinite2003}

As mentioned previously, there exist various formulations of the FVM that extend far beyond the original restriction of being first-order accurate.
The cell-averaged flux is then determined in terms of reconstructing polynomials that in theory can be of arbitrary order.
Such approaches per se do not fit well into the above given DG scheme but do however achieve similar results.

\section{Method for Assembling Numerical Schemes}

The overarching goal of this section is to identify a suitable combination of numerical schemes for a given multiphysics problem that is stable and accurate on the one hand, but also performant with regards to a specific choice of hardware on the other hand.

With the set relations between methods discussed in section \ref{sec:baseline}, we can now use the simplifications and thus computational advantages that each scheme presents.
That is, we follow the guideline to impose as many restrictions as possible whilst sustaining enough degrees of freedom to accurately capture the behavior of a given PDE.
In this way, we aim to provide the application expert with a recommendation on which schemes to use for a particular computational problem.
Our key concerns are, above all, to make this recommendation unambiguous.
Also, we focus on providing recommendations that will produce a stable solution and do not require tuning of artificial parameters.
If at least either of both requirements were not given, the usefulness of this method would be lost as one would have to undergo substantial experimentation to attain a valid solution.
Thus, by proposing such a method, we aim to give a sensible tradeoff between practicality and reproducibility on the one hand and utmost performance at the cost of possibly many model-building iterations and tuning on the other hand.

\subsection{Preliminary Assumptions}\label{sec:restrictions}

As a starting point, it has to be stated that encompassing the entire state of research on such schemes would be an impossible task.
The likewise formalization of a common framework is equally challenging as a consequence and thus not considered in this work.

Instead, we follow the path of introducing some restrictions that are on the one hand enough to construct a unifying scheme but on the other hand not too strict such that the efficient solution of real-world problems would be out of scope.

Thus, we propose the following restrictions to arrive at a one-to-one choice of numerical schemes:

\begin{enumerate}
    \item Only Bubnov Galerkin schemes are considered, that is, we omit Petrov Galerkin methods.
    The former restricts the choice of test space to be identical to the trial space.
    As such, we omit schemes that for instance use weighted functions or stencils to account for flow fields.
    An example of such schemes would be the Streamline Upwind Petrov Galerkin (SUPG) method \cite{brooksStreamlineUpwindPetrovGalerkin1982}.
    This restriction is essential to obtain an unambiguous choice of method, as the notion of Petrov Galerkin methods does not imply any particular choice of function space.
    \item We omit function spaces for approximation other than the $L^2$ and $H^1$ Sobolev spaces.
    There exists a vast variety of so-called Mixed Finite Element schemes that use Finite Elements based on different or composite function spaces with unique properties \cite{ernTheoryPracticeFinite2004}.
    For example, one may construct function spaces that can exactly fulfill divergence-free properties ($H(\text{div})$) or conditions based on the rotation of a field ($H(\text{curl})$).
    The specific choice of Finite Element then would require a considerable amount of expertise and would warrant a complex decision process of its own.
    We thus focus on scalar-, vector- and tensor-valued Lagrange elements solely.
    They have been shown to encompass a similar solution space as well and perform comparably for fluid and electromagnetic problems \cite{cockburnLocallyDivergencefreeDiscontinuous2004,hughesStabilizedMixedDiscontinuous2006}.
    \item Closely related to the previous statement, we restrict the solution space further by requiring that only Finite Elements utilizing Lagrange polynomials should be used.
    As the standard scalar- and vector-valued $\mathbb{P}^k$ and $\mathbb{Q}^k$ Finite Elements, being by far the most popular choices use exactly this family of polynomials, this requirement is weaker in practice than it might seem at first glance \cite{lagrange-0,lagrange-1}.
    \item We impose a coarse taxonomy to classify the qualitative behavior of a given PDE, that is, we specify limits regarding the leading coefficients of the differential operators.
    This should indicate whether the physical process described by the PDE is either more dissipative or more convective by nature. 
    We thus introduce a more physical interpretation than the considerably stricter coercivity measures employed in functional analysis.
    Our taxonomy closely follows the classes that were proposed by Bitsadze \cite{bitsadze1988some}.
    We do not claim this classification to be universally accurate.
    In practice, it has been shown however that having discrete cut-off values to disambiguate classes of PDE eases the choice of numerical scheme for application experts considerably.
    Hence, we choose to follow this path despite some shortcomings regarding generality.
    \item We only investigate systems of PDEs with differential operators up to second order.
    These are most common within physical processes and enable a wider range of numerical schemes to be used.
    For instance, equations of higher order such as the Cahn-Hilliard equation, would require the use of Finite Elements where up to third-order derivatives are defined.
    Such elements of high continuity are cumbersome to derive and are rarely used.
    Instead, we propose that in such cases the system should be reformulated as a mixed problem, where in the mentioned example one could represent the quantity of interest as two fields with second derivatives each.
    This technique is also regularly used in practice.
    \item For Finite Volume methods, we use the cell-centered variant as already outlined in section \ref{sec:fvm} instead of the vertex-centered or cell-vertex formulation.
    This is due to this form being the most popular choice in practice.
    Furthermore, it resembles the other schemes more naturally as has been shown in previous sections.
\end{enumerate}

\subsection{PDE Classification}\label{sec:pde-classification}

To find a numerical scheme that produces stable results, knowing the qualitative behavior of the system oftentimes is a necessity.
In particular, this means that the specific capabilities that a chosen numerical scheme possesses need to reflect the properties that the system presents.

We illustrate this by example. We once again investigate the simple advection equation (Eq. \ref{eq:advection-strong}), which is known to be first-order hyperbolic.
Such systems are prone to either preserving or even amplifying discontinuities given in the initial condition and thus the capability of accurately representing these should be incorporated into the choice of numerical schemes.
Suitable candidates would then be a Finite Volume or Discontinuous Galerkin Method.
However, the Finite Difference Method using a centered stencil or the Continuous Galerkin Method would give suboptimal results. The strong imposition of continuity in the domain would then yield spurious oscillations that affect stability.

We hence require the system of PDEs to be classified firsthand. We follow the popular taxonomy of second-order PDEs that can for example be found in the book by Bitsadze \cite{bitsadze1988some}, but follow a more general method for determining the appropriate class \cite{pinchoverIntroductionPartialDifferential2005}. 
That is, we define a singular governing equation in the form of a PDE to be either elliptic, parabolic or hyperbolic, depending on the shape that its characteristic quadric takes in space.
More formally, given a differential operator $L$ of the form:
\begin{equation}\label{eq:pde-classification}
    L[u] = \sum\limits_{i,j=1}^{n} a_{ij} \frac{\partial^2 u}{\partial x_i \partial x_j}
\end{equation}
where $x_i$ are the dependent variables and $a_{ij}$ is the matrix forming the coefficients of the highest spatial derivatives.
Considering the eigenvalues $\lambda_i$ of $a_{ij}$, $L$ is called

\begin{itemize}
    \item elliptic, if all $\lambda_i$ are either positive or negative,
    \item parabolic, if at least one eigenvalue is zero and all others are either positive or negative,
    \item hyperbolic, if at least one eigenvalue is positive and at least one is negative.
\end{itemize}

The characterization of first-order differential operators is more straightforward, however.
It can be shown that first-order PDEs with constant, real coefficients are always hyperbolic.
This condition is met for most cases relevant to engineering or physical applications.
More precisely, a first-order PDE is hyperbolic, if the resulting Cauchy problem is uniquely solvable.
In the case of real, constant coefficients, the polynomial equation for each variable has to admit $n$ solutions for an equation of order $n$ while keeping all other variables constant.
In the present case, this is trivially true.

We apply this classification for each governing equation of the independent variables for a given multiphysics problem.
In practice, one may oftentimes identify the class by the differential operators that frequently appear in a given PDE.
For example, a PDE that only has a laplacian as a spatial differential operator - such as the Laplace equation $\Delta u = 0$ or the heat equation $\partial_t u - \Delta u = 0$ exhibits dissipative behavior and is prototypical for elliptic and parabolic PDEs.
Oftentimes, one can easily identify a differential operator as parabolic if it has an elliptic operator in its spatial derivatives and an additional temporal derivative, as is exactly the case for the heat equation.

Both the abovementioned classes of PDE are dissipative, the reason being that PDEs of second order can only have discontinuous derivatives along their characteristics.
Since elliptic differential operators lack any characteristics, they strictly admit smooth solutions in that sense \cite{pinchoverIntroductionPartialDifferential2005}.
Thus, we associate this qualitatively dissipative behavior with elliptic and parabolic PDEs as defined above.

\begin{figure}[htbp!]
    \centering
    \resizebox{\linewidth}{!}{%
\begin{tikzpicture}[node distance = 2cm]

\node (toplevel) [element] {PDEs of (up to) second order};
\node (continuous) [element, below of=toplevel, xshift=-4cm] {Diffusive / Continuous};
\node (discontinuous) [element, below of=toplevel, xshift=4cm] {Convective / Discontinuous};

\node (elliptic) [element, below of=continuous, xshift=-2cm] {Ellptic};
\node (parabolic) [element, below of=continuous, xshift=2cm] {Parabolic};
\node (1hyperbolic) [element, below of=discontinuous, xshift=-2cm] {First order hyperbolic};
\node (2hyperbolic) [element, below of=discontinuous, xshift=2cm] {Second order hyperbolic};

\node (character) [description,left of=continuous, xshift=-3.5cm] {Character of the PDE / solution};
\node (type) [description,below of=character] {Class of the PDE};

\draw[arrow] (toplevel) -| (continuous);
\draw[arrow] (toplevel) -| (discontinuous);
\draw[arrow] (continuous) -| (elliptic);
\draw[arrow] (continuous) -| (parabolic);
\draw[arrow] (discontinuous) -| (1hyperbolic);
\draw[arrow] (discontinuous) -| (2hyperbolic);

\begin{pgfonlayer}{background}
    \node [fill=orange!30,fit=(character) (discontinuous)] {};
    \node [fill=blue!30,fit=(type) (2hyperbolic)] {};
\end{pgfonlayer}

\end{tikzpicture}
}%end resizebox
    \caption{Classification of PDEs up to second order by qualitative nature and types following \cite{bitsadze1988some}.}
    \label{fig:pde-classification}
\end{figure}
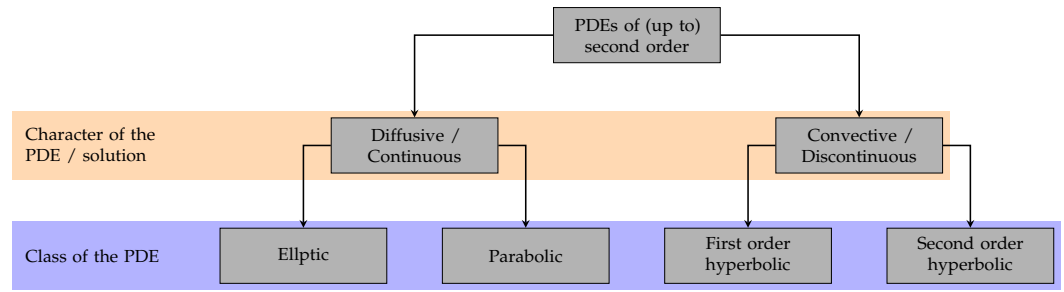

However, the advection equation (Eq. \ref{eq:advection-strong}) only has the gradient as a spatial differential operator, representing purely convective behavior.
Exactly this behavior of transporting information through the domain with finite speed is associated with the wave-like character of hyperbolic equations.
Figure \ref{fig:pde-classification} gives an overview of the classes of PDEs considered.

In alignment with the postulate at the beginning of this section, we aim to solve a given class of PDEs with as few degrees of freedom as possible whilst not over-constraining the solution.

Most importantly, discontinuities that might appear in the solution should be properly accounted for and reflect the choice of numerical scheme.
The direct consequence is that methods enforcing continuity should be used for problems that qualitatively exhibit high regularity and continuity.
From the previous discussion, it becomes apparent that this is the case for Finite Differences and Continuous Galerkin Finite Elements.
Problems that either conserve or even develop shocks however should be solved using methods that naturally allow for such.
This means that either Finite Volumes or Discontinuous Galerkin Finite Elements suit this requirement most naturally.

\subsection{Domain Geometry}\label{sec:domain-geometry}

As discussed in section \ref{sec:fdm}, the discretization using Finite Differences inherently assumes an even grid with uniform spacing between nodal points.
The direct consequence of this simplification is that assembly can be done in the computational domain directly and in an equal manner for every node point.

In general, if the domain has a particularly simple shape, for example, a hypercube and does not contain any holes, it can be triangulated using a cartesian grid.
Thus, if the discrete domain fulfills these conditions and the differential operators form an elliptic or parabolic PDE, using the FDM to efficiently assemble the global system is advisable.

For FVM, CGM and DGM, regularity of the computational domain, in general, does not pose any considerable advantages that may accelerate the assembly of the discretized system.

\subsection{PDE Linearity}\label{sec:pde-linearity}

Another crucial property to assess is its linearity.
In this case, we disambiguate strictly linear, semilinear, quasilinear and fully nonlinear equations, following the definition given in Evans \cite{evansPartialDifferentialEquations2010}:

A $k$-th order Partial Differential Equation of the form $$F(D^k u(x), D^{k-1} u(x),...,Du(x), u(x), x) = 0$$ is called:
\begin{enumerate}
    \item linear, if it has the form
    $$ \sum\limits_{\lvert \alpha \rvert \leq k} a_\alpha (x) D^\alpha u = f(x) $$
    for given functions $a_\alpha (\lvert \alpha \rvert \leq k), f$. The PDE is homogeneous if $f \equiv 0$. 
    \item semilinear, if it has the form
    $$ \sum\limits_{\lvert \alpha \rvert = k} a_\alpha (x) D^\alpha u + a_0 (D^{k-1}u,...,Du,u,x) = 0 $$
    \item quasilinear, if it has the form
    $$ \sum\limits_{\lvert \alpha \rvert = k} a_\alpha D^\alpha u (D^{k-1}u,...,Du,u,x) + a_0 (D^{k-1}u,...,Du,u,x) = 0 $$
    \item The PDE is fully nonlinear if it depends nonlinearly upon the highest-order derivatives.
\end{enumerate}
While linearity does not pose much of a problem for elliptic or parabolic equations, it plays an important role in whether a discretization is stable for hyperbolic equations.
The theory of nonlinear flux limiters is in general well researched for DG methods and largely profits from extensive developments that originally stem from the FVM.
However, accurate computation and implementation remain to be a hurdle in practice.
There have thus been several approaches to circumvent this issue, for example, by switching to an FV scheme in regions where there might be problems regarding the stability of the solution \cite{maltsevHybridDiscontinuousGalerkinfinite2023,sonntagShockCapturingDiscontinuous2014}.

As the overarching goal of this method is to provide straightforward guidance for end users, we will omit such approaches that must in most cases be implemented in a custom and rather particular fashion in favor of simplicity.

We thus recommend that, for equations where the solution is not likely to require many nonlinear iterations per time step, one may safely use a DG scheme.
In other cases where stability cannot be assured universally, one should rather switch to a Finite Volume formulation that may be overly diffusive, but on the upside is guaranteed to yield a stable solution.

\subsection{Computing Environment}\label{sec:hardware}

Within the last decade, advancement of computer hardware has been known to slowly hit the so-called \textit{memory wall} \cite{mckeeMemoryWall2011}.
That is, applications tend to be bound by the capability of the hardware to transfer memory instead of performing arithmetic operations.
This in particular holds for numerical simulations that are performed using many workers or large problems.
In such cases, the evaluation of sparse matrix-vector products poses high loads regarding memory bandwidth \cite{arndtExaDGHighOrderDiscontinuous2020}.

Then, there are numerical schemes that naturally lend themselves toward parallelism and others that are more memory-bound by design.
Thus, for a given computing hardware that puts enough emphasis on massive parallelism and two numerical schemes performing (nearly) identically, one should prefer the one that handles parallelism better.
We thus naturally arrive at the question of where one should disambiguate between massively parallel and other, regular hardware.

There are essentially two factors that would affect such a classification.
First, the hardware architecture itself plays an important role.
We may on the one hand solve a PDE on the classic CPU architecture that is capable of performing arithmetic on many precision levels and use many specialised instruction sets, such as AVX or fused multiply-add (FMA).
Another possibility is the use of highly parallel computing units, such as general-purpose graphics processing units.
Those however have a memory layout and instruction set that is much more tailored toward one purpose.
In the case of a GPU, this is medium to low precision operations with comparably low memory intensity but instead high arithmetic effort.

The other deciding factor is the amount of workers involved in the simulation process.
The more workers exist, the more processor boundaries are present and thus more information has to be shared between processors.
For some schemes, this overhead due to the exchange of memory between workers can become prohibitive.
Within the Finite Volume Method, for instance, parallel efficiency measured in GFlops/s starts to drop notably within the regime of 50 to 100 workers \cite{marshallFinitevolumeIncompressibleNavier1997}.
The quantitative drop-off also depends on the specific implementation, since other authors report slightly different results. 
Fringer et al. for instance note a decline in parallel efficiency for a Finite Volume solver starting at 32 workers \cite{fringerUnstructuredgridFinitevolumeNonhydrostatic2006}.
Thus, as a general guideline, we recommend employing methods that are suited for highly parallel environments at roughly 50 or more CPU workers.
For execution on massively parallel architectures, such as GPUs, the switch to such algorithms is considered necessary to obtain good efficiency.

\subsection{Problem Scale}\label{sec:problem-scale}

Another deciding factor for whether adaptivity is needed or not is the presence of multiple length scales in a multiphysics model.

We follow the definition given in \cite{eMultiscaleModeling2011,ePrinciplesMultiscaleModeling2011} and characterize a PDE-governed problem to have a Multiscale nature if models of multiple spatial or temporal scales are used to describe a system.
Oftentimes, this is the case if equations are used that originate from different branches of physics, such as continuum mechanics versus quantum mechanics or statistical thermodynamics.

This may on the one hand be a physical process with slow and fast dynamics, for example, in chemical reaction networks.
Then, the multiscale nature shows itself in the time domain of the problem.
Another example commonly encountered problem in alloy design is the evolution of the temperature field and phase kinetics during heating and solidification.
In this case, various length scales can be involved, such as in processes involving laser heating.
The temperature gradients then involve resolutions at a scale of around \SI{1e-5}{\metre}, whereas the width of a solidification front rather goes down to a sub-micrometer scale, that is, around \SI{1e-7}{\metre} \cite{zimbrodModellingMicrostructuresInsitu2021}.
About the previous definition, we have one model that is governed by laws of macroscale thermodynamics (that is,  the heat equation).
The other part of physics present is typically described by the evolution of a phase field.
The corresponding equations of this model are however derived from the formulation of a free energy functional from Landau theory \cite{landauTheoryPhaseTransitions1937}.

Due to the wide variety of physical processes and combinations thereof, formulating general criteria for the presence of a multiscale problem from a mathematical point of view is challenging.
To the knowledge of the authors, there do not exist any metrics in the literature that would enable such a classification.
We instead rely on the knowledge of the application expert who we assume to be familiar with the physics that should be captured.
For a rough disambiguation, however, one may use the definition given above.

Such multiscale phenomena are prohibitively expensive to resolve on a uniform mesh due to the nonnegligible difference in the dynamics of the system.
One option to efficiently resolve the physics at multiple scales is to employ different grids and solve the resulting problem in parallel.
This has, for example, been done for the above-mentioned case, specifically metal additive manufacturing \cite{ghanbariAdaptiveLocalglobalMultiscale2020}.

A rather effective, alternative approach is the modification of the governing equations such that they become tailored to a specific numerical scheme.
For instance, the well-known phase field model has been adapted using specialized stencils to the FDM such that spurious grid friction effects are eliminated \cite{fleckSharpPhasefieldModeling2023,fleckFrictionlessMotionDiffuse2022}. 
This approach however requires extensive knowledge about the numerics as well as the physical nature of a given problem.

Another possibility that requires fewer adaptions of the code to the specific problem is to make use of grid adaptive algorithms.
This approach for the problem presented is a popular alternative and has been implemented multiple times \cite{proellHighlyEfficientComputational2023,olleakEnablingPartScaleScanwise2022,olleakPartscaleFiniteElement2020,olleakScanwiseAdaptiveRemeshing2020a}.
Thus, grid adaptivity plays a key role in creating solutions to such problems, if the domains are not to be resolved on different discretizations entirely.
Numerical methods as a consequence need to reflect on this requirement and as such, Finite Difference methods are not suitable for such types of problems.

CG Finite Element methods do enable grid as well as polynomial degree adaptivity.
Yet, the imposition of hanging node constraints is oftentimes not trivial.
Though there have been considerable strides toward easy and intuitive handling of hanging nodes for continuous elements \cite{solinArbitrarylevelHangingNodes2008,bangerthDataStructuresRequirements2009}, these methods naturally fall short of the inherently decoupled nature of DoFs present in discontinuous methods.

Whereas grid adaptivity is easily realizable within FVM, there is little room for adaptivity regarding the order of approximation and can at best be achieved using varying reconstruction stencils \cite{shuHighorderFiniteDifference2003}.

By far, the most naturally suited method for h- as well as p-adaptivity is the DG FEM.
The locality of DoFs enables the splitting of cells without the need for hanging node constraints.
The same argument applies to altering the degree of a Finite Element, as additional DoFs within the cell need not be attached to a counterpart on its neighbors.

\subsection{Summary}

We may now condense the various aspects of choosing appropriate numerical schemes as follows into a unifying method, given the restrictions we posed in section \ref{sec:restrictions}.

First, we take three sets of inputs that are of practical relevance: the mathematically formulated, continuous problem, the computational domain that one wishes to solve the former on, and the configuration of the target hardware.

To design the intended decision process, we start by evaluating the decision metrics that impact the target scheme in the most general manner at first.
The general question of whether the prescribed system of PDEs requires an efficient solution on a large scale fulfills this requirement here.
By the term large scale, we understand state-of-the-art computing hardware on massively parallel architectures.
That decision in turn is influenced by two factors: One may directly intend to efficiently solve the system of PDEs on that hardware, or the multiscale nature of the problem demands such a computing environment.
If either is the case, solving the entire system using the HDG method is advisable due to the resulting locality of the problem.

The remaining parts of the decision process depend on the class of PDE present.
From here on, we operate in a field-wise manner and classify the system of PDEs for each independent variable separately.
If a PDE is convective in character, that is, hyperbolic, we recommend the use of numerical schemes that incorporate discontinuous approximations.
But, if a problem is diffusive by nature, the solution will be continuous and thus the use of continuous approximations is more advisable.

In the case of the former, following the discussion in section \ref{sec:pde-linearity}, a final disambiguation must be made regarding linearity.
If the PDE is linear or semilinear, a DG scheme can be applied due to the unlikeliness of stability issues.
Otherwise, the use of a simple FV scheme is more advisable to obtain a stable solution without having to iterate through many different choices of flux limiters in a trial-and-error fashion.

Regarding the continuous schemes, as has been explained in sections \ref{sec:fdm} and \ref{sec:domain-geometry}, the configuration of the domain geometry plays an important role in the efficiency of the overall scheme.
If the domain is cartesian, irrespective of dimensionality, the FDM can deliver accurate results with a considerably decreased amount of arithmetic operations.
The conceptual flexibility of the FEM regarding the domain is then unnecessary.
In the other case though where the domain is topologically more complex, relying on FEM algorithms that account for the necessary global mappings is more appropriate.
It would of course be possible to identify a middle ground between both schemes, for example, when a simple and prescribed transformation can be applied to the entire domain.
This would for example be the case for systems that can be described by polar coordinates.
However, few computer codes implement such functionality.
As the focus of this method lies on practicality and usefulness, we rather choose a method that can make use of widespread and established computer codes and thus omit these possibilities.

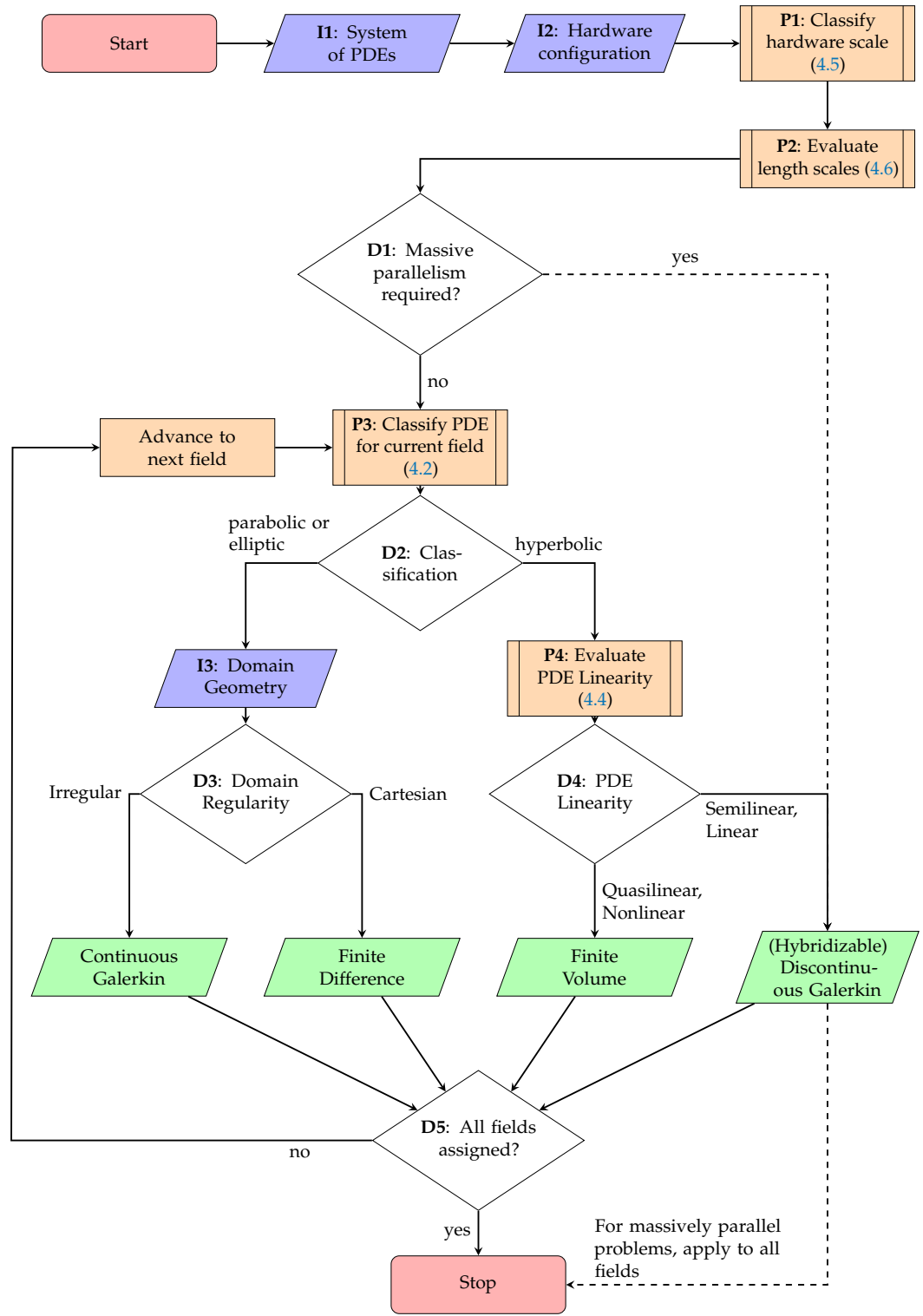
\begin{figure}[htbp!]
    \centering
    \resizebox{\linewidth}{!}{%
\begin{tikzpicture}[node distance = 2cm]

\node (start) [startstop] {Start};
\node (input-pde-system) [io, right of=start, xshift=2cm] {\textbf{I1}: System of PDEs};
\node (input-compute) [io, right of=input-pde-system, xshift=2cm] {\textbf{I2}: Hardware configuration};
\node (proc-compute) [subprocess, right of=input-compute, xshift=2cm] {\textbf{P1}: Classify hardware scale (\ref{sec:hardware})};
\node (proc-adaptivity) [subprocess, below of=proc-compute] {\textbf{P2}: Evaluate length scales (\ref{sec:problem-scale})};
\node (dec-compute) [decision, left of=proc-adaptivity,xshift=-5cm, yshift=-2cm] {\textbf{D1}: Massive parallelism required?};
\node (classify-pde) [subprocess, below of=dec-compute,yshift=-1cm] {\textbf{P3}: Classify PDE for current field (\ref{sec:pde-classification})};
\node (advance) [process, left of=classify-pde, xshift=-2cm] {Advance to next field};

\node (dec-class) [decision, below of=classify-pde] {\textbf{D2}: Classification};
\node (geometry) [io, below of=dec-class,xshift=-3cm] {\textbf{I3}: Domain Geometry};
\node (proc-linearity) [subprocess, right of=geometry, xshift=4cm] {\textbf{P4}: Evaluate PDE Linearity (\ref{sec:pde-linearity})};
\node (dec-hyperbolic) [decision, below of=proc-linearity] {\textbf{D4}: PDE Linearity};

\node (dec-geometry) [decision, below of=geometry] {\textbf{D3}: Domain Regularity};

\node (out-fdm) [res, below of=dec-geometry,xshift=2cm, yshift=-1cm] {Finite Difference};
\node (out-cgm) [res, below of=dec-geometry,xshift=-2cm, yshift=-1cm] {Continuous Galerkin};
\node (out-fvm) [res, below of=dec-hyperbolic, yshift=-1cm] {Finite\\Volume};
\node (out-dgm) [res, right of=out-fvm, xshift=2cm] {(Hybridizable) Discontinuous Galerkin};

\node (dec-done) [decision, below of=out-fdm,xshift=2cm,yshift=-1cm] {\textbf{D5}: All fields assigned?};
\node (stop) [startstop, below of=dec-done,yshift=-0.5cm] {Stop};

\draw[arrow] (start) -- (input-pde-system);
\draw[arrow] (input-pde-system) -- (input-compute);
\draw[arrow] (input-compute) -- (proc-compute);
\draw[arrow] (proc-compute) -- (proc-adaptivity);
\draw[arrow] (proc-adaptivity) -| (dec-compute);
\draw[dashed,arrow] (dec-compute) -| node[anchor=south,pos=.25] {yes} (out-dgm);
\draw[arrow] (dec-compute) -- node[anchor=west] {no} (classify-pde);
\draw[arrow] (classify-pde) -- (dec-class);
\draw[arrow] (dec-class) -| node[anchor=south,pos=.25, text width = 1.8cm] {parabolic or elliptic} (geometry);
\draw[arrow] (geometry) -- (dec-geometry);
\draw[arrow] (dec-geometry) -| node[anchor=east] {Irregular} (out-cgm);
\draw[arrow] (dec-geometry) -| node[anchor=west] {Cartesian} (out-fdm);
\draw[arrow] (dec-class) -| node[anchor=south,pos=.25] {hyperbolic} (proc-linearity);
\draw[arrow] (proc-linearity) -- (dec-hyperbolic);
\draw[arrow] (dec-hyperbolic) -- node[anchor=west, text width=2cm] {Quasilinear, Nonlinear} (out-fvm);
\draw[arrow] (dec-hyperbolic) -| node[anchor=north,pos=.25, text width = 2cm] {Semilinear, Linear} (out-dgm);
\draw[arrow] (out-cgm) -- (dec-done);
\draw[arrow] (out-fvm) -- (dec-done);
\draw[arrow] (out-fdm) -- (dec-done);
\draw[arrow] (out-dgm) -- (dec-done);
\draw[arrow] (dec-done)  -|  node[anchor=north,pos=.1] {no} ([shift={(-6.2cm,0mm)}]dec-done.west)-- ([shift={(-5.5cm,0mm)}]classify-pde.west)|- (advance);
\draw[arrow] (advance) -- (classify-pde);
\draw[arrow] (dec-done) -- node[anchor=east] {yes} (stop);
\draw[dashed,arrow] (out-dgm) |- node[anchor=south,text width=3.5cm,pos=.75] {For massively parallel problems, apply to all fields} (stop);

\end{tikzpicture}
}%end resizebox
    \caption{Graphic summary of the proposed process for choosing appropriate numerical schemes. Inputs (\textbf{I}) are given by purple trapezoids, decision points (\textbf{D}) by white diamonds and processes (\textbf{P}) by orange rectangles. Processes with additional vertical bars denote more complex processes and have references to their respective sections. Results are shown in green trapezoids.}
    \label{fig:decision-process}
\end{figure}

As a result, we obtain one process that guides the user through iteratively selecting the most appropriate combination of numerical schemes for a given, fixed and well-defined set of inputs.
This method may be summarised in a flow chart, which is depicted in Figure \ref{fig:decision-process}.

\section{Examples}

The purpose of this section is to walk through the proposed method employing two simple example PDEs.
Although these are not multiphysics problems, they may be combined in theory.

\subsection{Allen Cahn Equation}\label{Examples-Allen-Cahn-Equation}

First, we consider the following scalar PDE together with zero flux boundary conditions to be imposed at the four borders of a rectangular domain $\Omega: [0;L]\times[0;W]$:
\begin{align}
    \frac{1}{K} {\partial_t\phi} - \Delta \phi  = -\frac{2}{\xi^2}{\partial_\phi g(\phi)}-\frac{\mu_{0}}{3\gamma\xi}{\partial_\phi h(\phi)} \quad \forall \phi  \in \Omega, \label{eq:allen-cahn} \\
    \frac{\partial \phi}{\partial n}  = 0 \quad \forall x \in \partial \Omega \label{eq:ac-neumann},
\end{align}
This equation is called the Allen Cahn equation and describes the time-evolution of a scalar, non-conserved order-parameter field $\phi$, as is often called the phase field. 
The equation is commonly used in the modeling of self-organised microstructure evolution or complex pattern formation processes, as driven by local thermodynamics and/or mechanics. 
The phase field variable $\phi$ can be understood as a coloring function that locally indicates the presence or absence of a certain phase or a certain material state within a given microstructure. 
For instance, in modeling of microstructure evolution during solidification, $\phi = 1$ may denote the local presence of the solid and $\phi = 0$ may denote the local presence of the liquid phase \cite{fleckSharpPhasefieldModeling2023,fleckFrictionlessMotionDiffuse2022}.
If applied to the description of crack propagation, the order-parameter field $\phi$ is understood as the local material state, which can be either broken $\phi =1$ or not $\phi = 0$ \cite{PilipenkoFleckEmmerich2011,FleckPilipenSpatsch122010}.

The scalar quantities $K$, $\xi$, $\mu_0$ and $\gamma$ are model constants that determine the evolution of the scalar field $\phi$, and we adopt the notation of \cite{FleckFedermannPogorelov2018}. 
The polynomials $g$ and $h$ on the right-hand side of Eq. \ref{eq:allen-cahn} pose a nonlinearity to the equation. 
Their derivatives are given by $\partial_\phi g(\phi) = 2 \phi (1 - \phi) (1 - 2\phi)$ and $\partial_\phi h(\phi) = 6 \phi (1 - \phi)$. 
In the following, we will gather those polynomial terms in the joint potential term $f(\phi)= {2}{\partial_\phi g(\phi)}/{\xi^2}+{\mu_{0}}{\partial_\phi h(\phi)}/{3\gamma\xi}$.
Further details on the parametrization of the model are given in the appendix table \ref{tab:ac-parameters}.

Concerning the Allen Cahn Equation, we will consider two different scenarios, highlighting different aspects of the physics behind the equation. 
For each of these two scenarios, we formulate quantitative measures to be able to quantitatively question the accuracy of the numerical solutions.

In the first scenario, we consider the motion of a planar interface between two phases at different energy density levels. 
The low energy phase is expected to grow at the expense of the high energy phase, which induces a motion of the interface between them at a velocity proportional to the constant energy density difference $\mu_0$.
The scenario is realized as a quasi 1D problem $[0;L=100]\times[0;W=1]$, where the interface normal direction is pointing in the x-direction and the use of simple von Neumann boundary conditions with zero phase field fluxes at the borders of the rectangular domain is legitimate. 
The realization of this scenario with tilted interface orientations including the formulation of appropriate boundary conditions on the borders of the rectangular domain is discussed in detail in \cite{fleckFrictionlessMotionDiffuse2022}.
In this highly symmetric quasi-1D case, the scenario can be quantitatively evaluated using the existing analytic solution for the phase field:
\begin{align}
    \phi(x,0) = \phi_0(x) = \frac{1}{2} \left(1 - \tanh{\frac{x - \Tilde{x}(t)}{\xi}}\right) , \label{eq:ac-analytical}
\end{align}
where the time dependence of the central interface position $\Tilde{x}$ is given by, $\Tilde{x}(t) = x_0 + {M \mu_0 t}  /{\gamma}$, with the initial position at $x_0 = 20$.
To investigate the impact of arithmetic complexity on computational efficiency, we will seek an approximation of a real-valued function $\phi (x,t)$ in one spatial dimension on a $[0;100]$ grid with equispaced vertices.

We can now start applying the proposed methodology, as described above in section \ref{sec:baseline} and following Figure \ref{fig:decision-process}.
That is, we follow the path of the flowchart from top to bottom.
We first classify the hardware scale according to P1 in the figure.
The given hardware architecture, that is, ~an 8-core CPU system, falls well below the established recommendation for the threshold of partitioned problems which is at least 50 workers.
Therefore, there is no need from a hardware side for massive parallelism.

Next, we investigate the problem scale with process P2.
As the problem is governed by one scalar equation and there are no sub-models involved as defined in section \ref{sec:problem-scale}. 
Concerning the length scales, the presented system exhibits one extra physical length scale and that is the width $\xi$ of the diffuse interface.
This extra physical length scale originates from the nonlinearity of the Allen Cahn Equation and complements the other length scales such as the dimensions of the domain as well as the grid spacing, both being more natural in the numerical solution of PDEs. 
This poses the issue of numerical resolution of the systems length scales, that is, ~both the domain dimensions, as well as the width $\xi$ of the diffuse interface, need to be properly represented on the discrete numerical grid \cite{fleckSharpPhasefieldModeling2023,fleckFrictionlessMotionDiffuse2022}.
However, the fact that the problem is quasi-one-dimensional restricts the computational demands of the scenario.
We thus arrive at the first decision point D1, where we can negate the necessity for massive parallelism.

The next process step P3 involves classifying the problem at hand, following the definition given in section \ref{sec:pde-classification}.
As dependent variables, we encounter the time $t$ as well as the spatial components $x$ and $y$.
The coefficient matrix $a_{ij}$, summing up all leading coefficients of second derivatives then becomes for the 2D case:
\begin{equation}
    a_{ij}^{AC} = \begin{bmatrix}
        0 & 0 & 0 \\
        0 & -1 & 0 \\
        0 & 0 & -1
    \end{bmatrix}
\end{equation}
In this case, enumerating the eigenvalues $\lambda_i$ is trivial, since $a_{ij}^{AC}$ is a diagonal matrix, and we have $\lambda_0 = 0, \lambda_1 = -1, \lambda_2 = -1$.
We thus find that one eigenvalue is zero since the temporal derivative is only of first order and all other eigenvalues are of the same sign.
Therefore, Eq. \ref{eq:allen-cahn} is a second-order PDE of parabolic type and we can proceed in D2 with the left branch.

Moving on in the decision process, we would next classify the problem domain in D3 given input I3.
As we use an equispaced grid in 1D, the discretization is cartesian and thus solving the problem using Finite Difference would be the best choice.
As there are no other fields to classify according to decision point D5, we conclude the decision process.
Within the unified methodological framework, we implement both FD and CG schemes and the scenario is comparatively solved using both schemes.
This allows us to compare the schemes concerning numerical resolution capabilities and to investigate differences in the mutual arithmetic complexity and their impact on efficiency. 

Evaluating the CG method requires the re-formulation of Eq. \ref{eq:allen-cahn} in its weak form, though.
The finite dimensional weak statement is then: Find $\phi_h \in V_h$, such that
\begin{equation}\label{eq:ac-weakform}
    \int_\Omega \frac{1}{M} v_h \partial_t \phi_h \,dx + \int_\Omega \nabla v_h \nabla \phi_h \,dx - \int_\Omega v_h f(\phi_h) \,dx = 0 \quad \forall v_h \in V_h
\end{equation}
Where we have already assumed the solution and test function to lie in the finite-dimensional subspace $V_h$.

Eq. \ref{eq:allen-cahn} requires the discretization of the Laplacian as its only differential operator.
The temporal derivative will be treated using the Method of Lines approach, that is,  we solve a large system of spatially discretized ordinary differential equations.

The Finite Difference discretization of the laplacian results in the well-known second-order central difference stencil:
\begin{equation} \label{eq:ac-stencil}
    \Delta \phi \approx \frac{\phi_{i+1} - 2 \phi_{i} + \phi_{i-1}}{\Delta x^2}
\end{equation}
The nonlinear right-hand side $f(\phi)$ must be updated every time step using the current value of $\phi$.
As such, we do not need to perform any assembly and can even avoid forming a global system of equations.
Instead, we rely on \ref{eq:ac-stencil} for the Laplacian, which can be handily vectorized.
There is also no need to perform any mapping between the reference and physical domain as explained in section \ref{sec:fdm}.

For the Finite Element discretization, we need to perform all these steps, resulting in a global nonlinear system of equations for each time step. The discrete form of Eq. \ref{eq:ac-weakform} then reads:
\begin{equation}
    \underline{\underline{M}} \partial_t \underline{\phi} + \underline{\underline{K}} \, \underline{\phi} - F(\underline{\phi}) = 0
\end{equation}
Where we introduced the mass matrix $\underline{\underline{M}}$ and the stiffness matrix $\underline{\underline{K}}$ for the Laplacian.
These represent the spatially discretized differential operators that act on the vector of degrees of freedom $\underline{\phi}$.
The algebraic terms that are nonlinear in $\underline{\phi}$ are gathered in the discrete vector $F(\underline{\phi})$.
For the sake of comparison regarding efficiency, we require the resulting fields of both schemes to be (nearly) identical apart from floating point errors. 

Given this requirement, we note that the Finite Difference formulation lacks an analogous term to the Finite Element mass matrix. 
We consequently require $M$ to be the identity matrix in an equivalent Finite Element formulation, given all other terms are equal.
The latter can easily be verified for a stiffness matrix assembled with first-order Lagrange polynomials and a collocation method.
The derivation of such an equivalent scheme has been covered in section \ref{sec:fdm}.
Using collocated Finite Elements is chosen here for the sake of comparison as well as for computational efficiency.
The resulting mass matrix can then be inverted trivially by taking the element-wise inverse instead of computing the full inverse.
Such an operation is considerably more expensive and should thus be avoided if possible.

To compare both schemes regarding efficiency, we implement both schemes from scratch within the Julia programming language \cite{bezansonJuliaFreshApproach2017}.
Due to its flexibility, high-level syntax and simultaneous, granular control over various performance aspects via its rich type system, Julia has gained considerable momentum in the past few years within the scientific community.
We carefully set up both schemes using analogous data structures to enable a side-to-side comparison of the computational complexity.
The most high-level parts of the codes are given in Listing \ref{lst:ac-code}.

\begin{figure}[!htbp]
\begin{minipage}[t]{0.5\textwidth}
\centering
\begin{minted}[breaklines,escapeinside=||,mathescape=true, linenos, numbersep=3pt, gobble=2, frame=lines, fontsize=\footnotesize, framesep=2mm]{julia}
struct FETriangulation{V,C}<:Triangulation
    vertices::V
    connectivity::C
    dim::Int
end

struct FiniteElement{E<:ElementType,
        P<:Primitive,B<:AbstractMatrix,
        Q<:AbstractVector,G<:AbstractArray}
    primitive::P
    element_type::E
    order::Int
    ndofs::Int
    basis_coeffs::B
    quadrature_nodes::B
    quadrature_weights::Q
    basis_at_quad::B
    grad_basis_at_quad::G
    grad_monomial_basis::G
end

struct AssemblyCache{T<:AbstractVecOrMat}
    coeffs::T
    loc::T
    glob::T
end

struct CGProblem{T<:Triangulation,J,P,C
      E<:FiniteElement,M<:AbstractMatrix,
      V<:AbstractVector,F<:Function}
    triangulation::T
    referenceElements::E
    detJ::J
    bilinearForm::M
    linearForm::V
    u::V
    massMatrix::Union{M,Nothing}
    parameters::P
    rhs::F
    cache::C
end

function (a::CGProblem)(du,u,p,t)
    mesh        = a.triangulation,
    element     = a.referenceElements
    K           = a.bilinearForm
    F           = a.linearForm
    detJ        = a.detJ
    params,     = a.parameters
    cache,rhs   = a.cache, a.rhs
    M           = a.massMatrix

    mul!(du,K,u,-1.,0.)
    # f(u) changes, thus we reassemble
    assemble_F!(F,cache,u,element,
            mesh,detJ,params,rhs,M)
    du .+= F

end
\end{minted}
\end{minipage}
\begin{minipage}[t]{0.5\textwidth}
\centering
\begin{minted}[breaklines,escapeinside=||,mathescape=true, linenos, numbersep=3pt, gobble=2, frame=lines, fontsize=\footnotesize, framesep=2mm]{julia}
struct FDTriangulation{V,D}<:Triangulation
    vertices::D
    h::V    # dx in each dim
    dim::Int
end






















struct FDProblem{T<:FDTriangulation,
        B<:Function,L<:Function
        V<:AbstractArray,P,BC<:Function}
    triangulation::T
    order::Int
    bilinearForm::B
    linearForm::L
    u::V
    parameters::P
    boundaryCondition::BC
end




function (a::FDProblem)(du,u,p,t)
    apply_bilinear! = a.bilinearForm
    apply_linear!   = a.linearForm
    apply_bc!       = a.boundaryCondition
    h               = a.grid.h
    params          = a.parameters




    a.bilinearForm(du,u,h)



    a.linearForm(du,u,a.parameters)
    a.boundaryCondition(du)
end
\end{minted}
\end{minipage}

\vspace{0.2cm}
\captionof{listing}{Top-level overview of the necessary data structures and the functions to update the semi-discrete systems for the CG FEM (left) and FDM (right).
In the case of the FDM, one can avoid assembling a global linear system entirely, thus the top-level data structure only holds the solution and grid as large arrays.
For the FEM, assembly on general grids in a matrix-free manner is far from trivial.
Additionally, the triangulation data structure is more complex due to the necessary topological information.
Furthermore, the reference FE needs to be stored and correctly mapped using Jacobian values. The full code is available in the code repository mentioned at the end of this article.}
\label{lst:ac-code}
\end{figure}

We also include the functions that are called within each time step to solve the semidiscrete system, to give a high-level view of which steps are necessary and how they are implemented in particular. 
Both semidiscrete systems use in-place operations to avoid memory allocations.
For the CG-FEM code, we implement a full mesh topology to solve the problem with a first-order method although both discretizations consist of cartesian meshes.
One could in this case assume a globally constant jacobian and thus save a considerable amount of arithmetic complexity.
However, this would skew the results regarding performance and would not make full use of the flexibility of the FEM.

It becomes immediately apparent from the comparison that solving the Allen Cahn equation using Finite Elements requires an assembly process that is noticeably more complex.
The only arrays that need to be stored for the FD version are the grid coordinates and the solution array.
Because the latter can be arranged in memory such that it represents the cartesian topology of the grid, one can simply point to the neighbors of a vertex in memory without having to look up the vertex-vertex connectivity.
This is not the case for the FEM.
Instead, we encounter an additional indirection through a cell-vertex list, where we gather all DoFs associated with the currently visited cell.

We furthermore cannot construct the global linear system at once, but need to go through the cell-wise assembly process which effectively leads to most of the non-zero matrix entries being visited multiple times.
This is in sharp contrast to the FDM where the global system is only present implicitly through functions that apply the laplacian stencil.
As a consequence, memory requirements are greatly reduced.

For transient problems, one needs to additionally make a suitable choice for the temporal discretization, that is,  the choice of method as well as the time step.
Here, we make use of the well-optimized Julia library \texttt{DifferentialEquations.jl} \cite{rackauckasDifferentialequationsJlPerformant2017}.
As an exemplary implementation of modern, high-performance codes for the solution of ordinary differential equations (ODE), this package offers various algorithms that are capable of adaptive time stepping such that an application expert does not need to provide any input regarding temporal discretization.
Here in particular, we can even make use of built-in heuristics that automatically select a suitable integration scheme, based on the supplied ODE problem \cite{rackauckas2019confederated}.
The resulting effort for the end user can be condensed to selecting a suitable numerical scheme for the spatial discretization as outlined by Figure \ref{fig:decision-process} and leave the problem of tuning the spatial discretization aside entirely.
For this particular problem, we prescribe the use of an adaptive, implicit, 4th-order Rosenbrock method for the temporal evolution of both FD and CG-FE systems to achieve a fair comparison between both solutions.
This solver is stable and third-order accurate when used on nonlinear parabolic problems \cite{rackauckasDifferentialequationsJlPerformant2017}.

\begin{figure}
     \centering
     \begin{subfigure}[t]{0.49\textwidth}
         \centering
         \includegraphics[width=\textwidth]{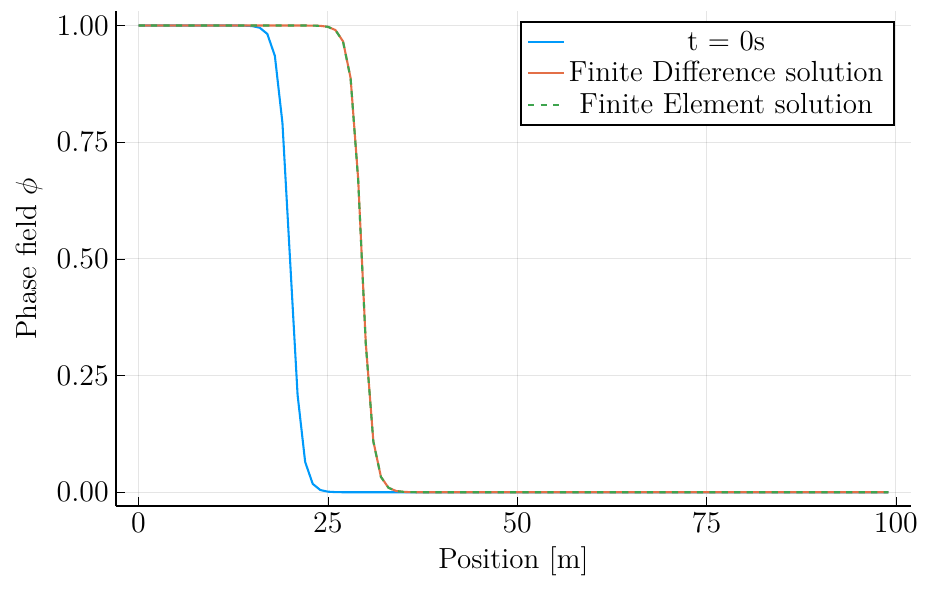}
         \caption{Evolution of the phase front at $t=100s$ with respect to the initial condition.}
         \label{fig:ac-phasefield}
     \end{subfigure}
     \hfill
     \begin{subfigure}[t]{0.49\textwidth}
         \centering
         \includegraphics[width=\textwidth]{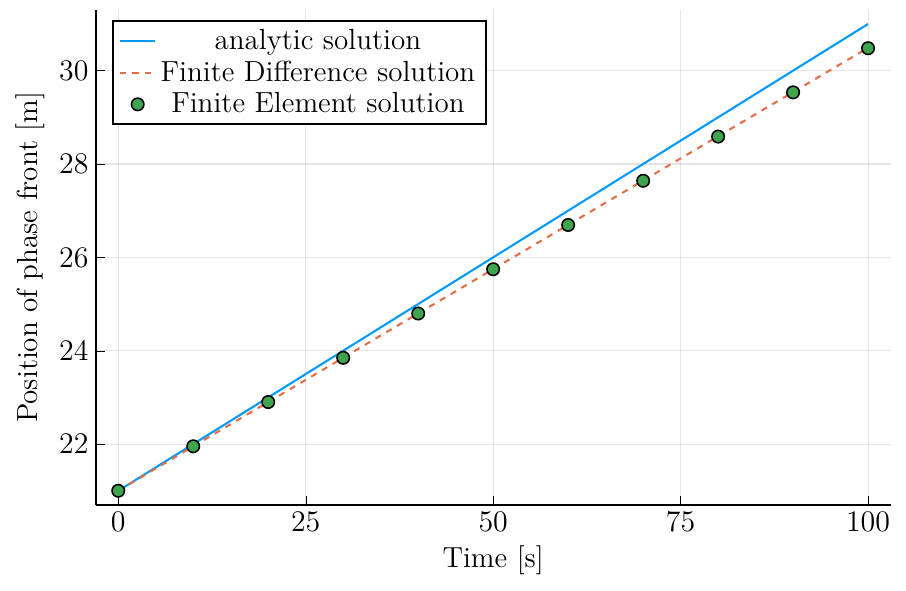}
         \caption{Comparison of Finite Difference and Finite Element solutions to the analytical solution given by Eq. \ref{eq:ac-analytical}}
         \label{fig:ac-positions}
     \end{subfigure}
    \caption{Solution of the one dimensional Allen Cahn equation using the Finite Difference and Finite Element method.}
    \label{fig:ac-solutions}
\end{figure}

Before comparing both schemes regarding computational efficiency, we first verify that the FD and collocated CG-FE schemes produce identical results. Figure \ref{fig:ac-solutions} shows the solutions of both schemes for solving the phase field evolution (Figure \ref{fig:ac-phasefield}) and for modeling interface position over time (Figure~\ref{fig:ac-positions}).

As can be observed, both schemes produce visually identical results.
The quantitative differences in the numerical results are minimal and can be attributed to floating point errors that accumulate over the process of time integration.
However, there is a considerable difference between the analytical and the numerical interface velocity, as visible in Figure \ref{fig:ac-positions}.
The reason for this discrepancy is grid friction, which results from the limited numerical resolution of the diffuse interface profile and could be reduced by increasing the dimensionless ratio $\xi /\Delta x = 1.5$, where $\Delta x$ denotes the grid spacing. \cite{fleckSharpPhasefieldModeling2023,fleckFrictionlessMotionDiffuse2022}. 
It is quite interesting to note, that the grid friction effect turns out to be so very similar for the two different numerical schemes in this case.

\begin{table}
\centering
\caption{Run times of the Finite Element and Finite Difference model of the 1D Allen Cahn equation. Both models were run on identical hardware and Julia 1.8.5 with LLVM 13.0.1 \cite{bezansonJuliaFreshApproach2017}. Time stepping was performed using the \texttt{DifferentialEquations.jl} library \cite{rackauckasDifferentialequationsJlPerformant2017}. Linear Algebra operations are performed using OpenBLAS \cite{wangAUGEMAutomaticallyGenerate2013,xianyiModeldriven} on a single-threaded Apple M1Pro ARM processor. Fast evaluation of fused array expressions is provided by the \texttt{Tullio.jl} library \cite{abbottMcabbottTullioJl2022}. Allocated Memory refers to the physical size of the problem-specific data structures given in \ref{lst:ac-code}. The sample size for each scheme is $n=100$.}
\label{tab:ac-performance}
\begin{tabular}{llll} 
    \toprule
    & \textbf{FDM} & \textbf{FEM}  & \textbf{Relative}  \\
    \midrule
Median run time             & 0.450 ms           & 9.503 ms             & 21.1x    \\
Mean run time $\pm 1\sigma$ & 0.578 ms ± 0.759 ms & 9.643 ms ±  0.741 ms & 16.7x    \\
Allocated Memory            & 1.446 kB            & 15.598 kB             & 10.8x  \\
\bottomrule
\end{tabular}
\end{table}

Furthermore, we point out that computational resource usage differs considerably. 
Table \ref{tab:ac-performance} reports some descriptive statistics on the performance of both implementations.
These differences in run times as well as memory consumption can be attributed to multiple factors.
First, the nonlinear right-hand side changes each time step and thus assembly has to be performed dynamically for the Finite Element Method. The Finite Difference Method in contrast can simply rely on point-wise evaluation of the strong form instead of numerically computing the weak form integrals.
Secondly, the Finite Difference Method does not need to perform any mapping during the time step as no assembly is required. During computation of the right-hand side integral, this is a necessity for the Finite Element Method.

The largest discrepancy however can be attributed to the fact that the Finite Difference Method can operate in a matrix-free manner due to the cartesian grid it is applied on. As all vertices are equispaced, there exists one global stencil that can be applied on each vertex independent of all other members of the grid. The Finite Element Method in contrast uses the grid topology to accumulate the weak form integrals into corresponding entries of the global system matrices and vectors. Thus, it always produces a typically very sparse global system that cannot be vectorized similarly.
It should be noted that the discrepancy in results should not be expected to be as drastic as shown for linear problems, as then the FEM does not require the re-assembly of the right-hand side.
The computational advantage then reduces to the matrix-free evaluation of the linear system.

\begin{figure}
    \centering
    \begin{subfigure}[t]{0.49\textwidth}
        \centering
        \includegraphics[width=\textwidth]{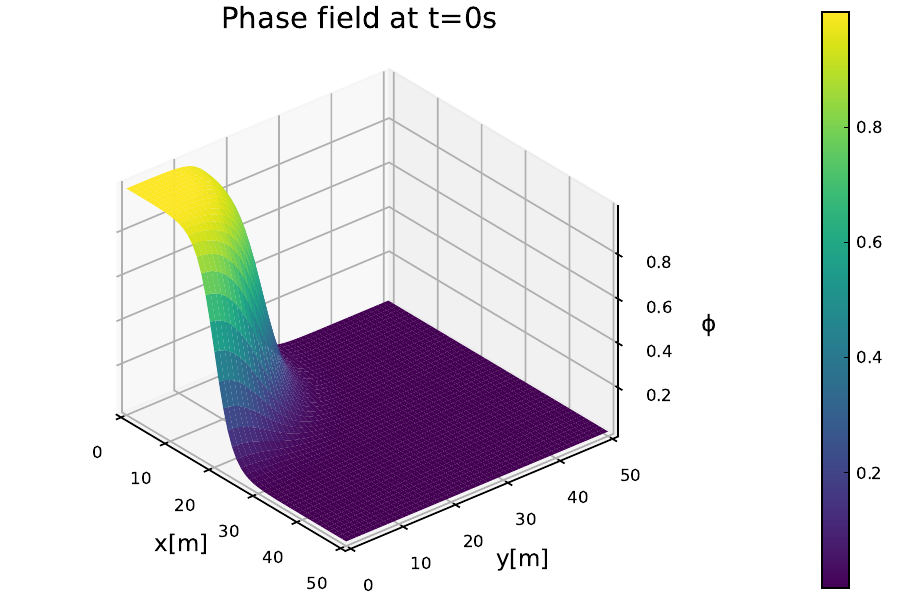}
        \caption{Initial configuration of the phase field. The domain shows a quarter slice of the nucleus.}
        \label{fig:nucleus-ic}
    \end{subfigure}
    \hfill
    \begin{subfigure}[t]{0.49\textwidth}
        \centering
        \includegraphics[width=\textwidth]{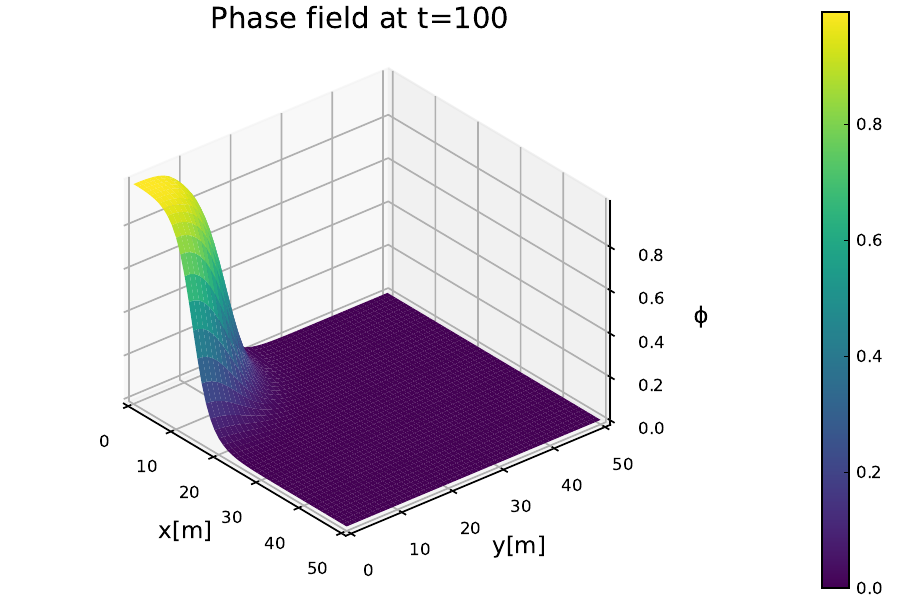}
        \caption{Phase field distribution within the quarter domain after 100 seconds. The nucleus has shrunk to a smaller radius whilst retaining the interface profile.}
        \label{fig:nucleus-final}
    \end{subfigure}
   \caption{Solution of the two-dimensional Allen Cahn equation using the Finite Difference and Finite Element method.}
   \label{fig:ac-nucleus}
\end{figure}

In the second scenario, we investigate a more practically relevant benchmark in two dimensions and turn to the well-known vanishing grain problem, leaving all other aspects of the problem as is.
Here, the dissolution of a circular-shaped nucleus under the interface energy density pressure under two-phase equilibrium condition $\mu_0=0$ is simulated. 
These dynamics are also governed by the Allen Cahn equation and denote the complementary physical effects as compared to the above scenario.
In a sharp interface picture, with a constant and isotropic interface energy density $\gamma$, we expect the temporal evolution of the grain radius to be given by:
\begin{equation}
    r(t) = \sqrt{R_0^2 - 2 M t}
\end{equation}
Where $R_0$ indicates the initial radius and $M$ is the phase field mobility.
Snapshots of the phase field $\phi$ at initial and terminal times are given in Figure \ref{fig:ac-nucleus}.

\begin{figure}
    \centering
    \includegraphics[width=0.7\textwidth]{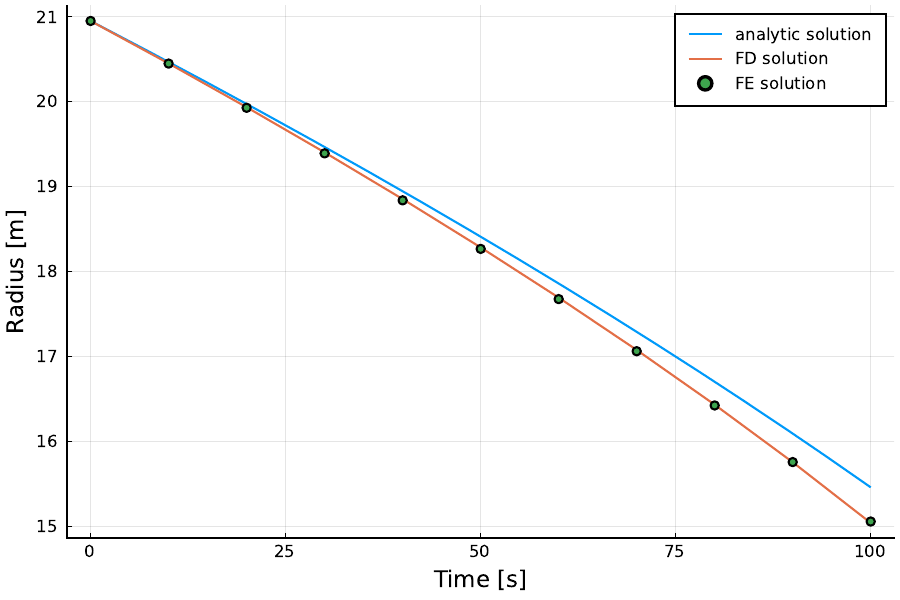}
    \caption{Evolution of radius over time of the vanishing grain problem. Both Finite Difference and Finite Element solutions show a considerable, accumulating error toward the analytical solution.}
    \label{fig:ac-radius}
\end{figure}

We also report the temporal evolution of the radius function for solving this scenario using both numerical methods in Figure \ref{fig:ac-radius}.
As in the one-dimensional simulation, the collocated FE and FD solutions behave identically to each other.
Both exhibit a notable discrepancy towards the sharp interface behavior, which again relates to known issues of finite numerical resolution in the phase field simulation \cite{fleckFrictionlessMotionDiffuse2022}.
We note at this point that the solution generated by a common FE model, i.e. with fully accurate quadrature as explained in section \ref{sec:fdm} also produces very similar results given all other parameters are chosen the same, albeit with a large computational disadvantage due to the full inversion of the resulting mass matrix.
The linked code repository at the end of this article contains the necessary data structures to reproduce these results.

\begin{table}
    \centering
    \caption{Run times of the Finite Element and Finite Difference model of the 2D Allen Cahn equation. Both models were run on identical hardware and Julia 1.8.5 with LLVM 13.0.1 \cite{bezansonJuliaFreshApproach2017}. Time stepping was performed using the \texttt{DifferentialEquations.jl} library \cite{rackauckasDifferentialequationsJlPerformant2017}. Linear Algebra operations are performed using OpenBLAS \cite{wangAUGEMAutomaticallyGenerate2013,xianyiModeldriven} on a single-threaded Apple M1Pro ARM processor. Fast evaluation of fused array expressions is provided by the \texttt{Tullio.jl} library \cite{abbottMcabbottTullioJl2022}. Allocated Memory refers to the physical size of the problem-specific data structures given in \ref{lst:ac-code}. The sample size for each scheme is $n=100$.}
    \label{tab:ac2d-performance}
    \begin{tabular}{llll} 
    \toprule
                                & \textbf{FDM}                 & \textbf{FEM}                   & \textbf{Relative}  \\
                                \midrule
    Median run time             & 4.743 ms            & 192.837 ms             & 40.7x    \\
    Mean run time $\pm 1\sigma$ & 4.731 ms ± 0.059 ms & 192.787 ms ±  0.219 ms & 40.7x    \\
    Allocated Memory            & 21.490 kB            & 494.730 kB              & 23.0x  \\
    \bottomrule
    \end{tabular}
    \end{table}

To assess the performance gap of both schemes for higher dimensions, we again benchmark both codes against each other.
The results are given in Table \ref{tab:ac2d-performance}.
Comparing the results from the 2D simulation benchmark in table \ref{tab:ac2d-performance} with its 1D counterpart (table \ref{tab:ac-performance}), we find that the discrepancy in performance becomes noticeably more drastic with increasing dimensionality of the problem.
This can be attributed to the increased scattering of DoFs in memory.
Thus, memory access is less stridden, increasing the lookup time.
For a hardware architecture that demands more parallelism and has a shared memory architecture, this could quickly evolve into a serious bottleneck.

\subsection{Two-phase Advection}\label{sec:advection-equation}

As a second model problem, we will investigate the advection equation in two dimensions.
This problem is well-studied in the literature and is known as challenging to solve accurately.
Due to the absence of dissipative terms, numerical algorithms oftentimes struggle to converge towards the entropy solution and either produce spurious oscillations, rendering the solution unstable or yield overly diffusive approximations, where conservation laws are violated \cite{levequeNumericalMethodsConservation1992}.
We choose this problem in particular due to being simple yet challenging enough to study.
In addition, the advection equation frequently arises in modeling multiple phases in an Eulerian framework and the motion of immersed immiscible fluids in general.
It is thus of high relevance in a multitude of multiphysics problems.

In particular, we investigate a pure advection problem involving two phases. 
We choose to describe the motion of two fluids and track the volume fractions $\alpha_i$, as is common for the Volume-of-Fluid (VoF) formulation:
\begin{align}
    \partial_t \alpha + \underline{u} \cdot \nabla \alpha = 0  \label{eq:vof-advection}\\
    \alpha_1 + \alpha_2 = 1 \label{eq:vof-sum} \\
    \Omega \in \left[0;5\right] \times \left[0;5\right] \\
    t \in \left[0;5\right]
\end{align}
The initial condition to this problem is given as a rectangle function that is one in the interval $x \in \left[2;3\right] \times \left[2;3\right]$ and zero everywhere else.
One may alternatively track only the motion of the interface using a coloring function $\phi$. 
This is common for the level set method, the governing equation however is the same as Eq. \ref{eq:vof-advection}.

We would like to solve this problem on three different architectures to showcase the effect of parallelism on the efficiency of numerical schemes.
The choices of hardware along with important quantities are given in Table \ref{tab:advection-hardware}.
\begin{table}[htbp!]
    \centering
    \caption{Hardware configurations for the advection equation model problem. The three setups mimic popular computing environments in applied settings: A mobile computer, a stationary workstation grade tower and a server tailored to numerical computing.}
    \begin{tabularx}{\linewidth}{LLLLLL} 
    \toprule
    \textbf{CPU Name} & \textbf{Number of Cores} & \textbf{Core Clock Speed [GHz]} & \textbf{Memory Size [GB]} & \textbf{Memory Bandwidth [GB/s]} & \textbf{Memory Speed [GHz]}  \\
    \midrule
    Apple M1 Pro      & 8 & 3.2 & 16 & 200 & 6.4 \\
    Intel Xeon W-2295 & 18 & 3.0 & 128 & 94 & 2.9 \\
    2x AMD EPYC 7763  & 128 & 2.45 & 512 & 204 & 3.2 \\
    \bottomrule
    \end{tabularx}
    \label{tab:advection-hardware}
    \end{table}
We once again follow the process summarized in Figure \ref{fig:decision-process}.
Regarding the system of PDEs (I1), Eq. \ref{eq:vof-sum} is simply an algebraic constraint and thus can be calculated in a simple postprocessing step.
Thus, \ref{eq:vof-sum} is not a governing equation in the sense of a PDE and $\alpha_2$ will consequently not be considered an independent variable, as detailed in section \ref{sec:pde-classification}.
We are then left to solve a single scalar advection equation for $\alpha_1$.

Proceeding in the flow chart, we next classify the hardware scales within P1.
Here we find that the last hardware configuration listed in Table \ref{tab:advection-hardware} necessitates the use of schemes that are tailored for high parallelism, as the given amount of 128 processes is above the specified regime where the use of parallelizable algorithms is worthwhile using.
Therefore, this configuration should be run using a Discontinuous Galerkin Method.
For both other configurations using 8 and 18 processes, this does not apply.
Continuing with process P2, we find that the given problem only exhibits one length scale and thus this criterion for parallelism can be omitted.
Thus, we arrive at decision D1 and find that for the problem statement involving the largest of the three computing architectures, the use of the Discontinuous Galerkin Method is advised.

For the remaining two configurations, we can proceed by classifying the PDE according to process P3.
With the temporal derivative and gradient as the only differential operators, Eq. \ref{eq:advection-strong} is a first-order PDE.
The advection velocity vector $\underline{u}$ has constant and real components.
Thus, following section \ref{sec:pde-classification}, we find that the advection equation presented here is hyperbolic and proceed with the right branch of the flow chart after decision D2.

Consequently, we need to evaluate the linearity of Eq. \ref{eq:vof-advection} for process P4 as a next step.
As the terms including the differential operators are linear and there is no right-hand side, it may be classified straightforwardly to be linear.
Due to its linearity, the most efficient choice for the remaining configurations turns out to be the Discontinuous Galerkin method as well.
As stated previously, the original two-equation system only consists of one PDE, and thus we conclude the decision process here, as all fields governed by a PDE have been assigned (decision D5).

In the following, we will compare this particular choice of method with the Finite Volume Method, which would be the next alternative and is in principle also well suited to tackle such problems.
The Continuous Galerkin and Finite Difference Methods do not lend themselves well to solving such equations and will thus be omitted from this benchmark.
In particular, the CG method is known to be unstable for first-order hyperbolic equations, as the stability of the scheme can be shown to be dependent on mesh size \cite{ernTheoryPracticeFinite2004}.

One must add that in principle, the Finite Difference Method could be applied here, where however two different limitations apply.
First, the only choice of stencil that would be stable for this equation is the forward difference (or upwind) approximation.
This choice however is not covered by the proposed decision process, as it is formally equivalent to a Continuous Petrov Galerkin Method.
One can show that this scheme corresponds to a simplified Streamline Upwind Petrov Galerkin (SUPG) method. As we have restricted ourselves for the sake of decidability to Bubnov Galerkin methods, this stencil is not admissible here.
Secondly, using the Finite Difference Method here implies the strict use of a cartesian grid.

We thus proceed to write the weak form for Eq. \ref{eq:vof-advection} by multiplying with a discrete test function $v_h$, integrating over the whole domain $\Omega$ and subsequently performing integration by parts. Find $\alpha_h \in V_h$ such that:
\begin{align}
    \int_\Omega v_h \frac{\partial \alpha_h}{\partial t} \,dx + \int_\Gamma v_h \alpha_h (\underline{u}\cdot\underline{n}) \,dS - \int_\Omega (\nabla v_h \cdot \underline{u}) \alpha_h \,dx = 0     \qquad \forall v_h \in V_h
    \label{eq:advection-weakform} \\
    \alpha_{h,\tilde{\Gamma}_{b,l}} = \alpha_{h,\tilde{\Gamma}_{t,r}} \quad \forall x \in \partial\Omega
\end{align}
Where $\Gamma$ denotes the union of all interior and exterior facets of the domain and $\tilde{\Gamma}$ are the subsets of the domain boundary $\partial \Omega$. In this case specifically, $\tilde{\Gamma}_{b,l}$ are the slave facets at the bottom and left boundary that the values of the slave facets from the top and right master facets $\tilde{\Gamma}_{t,r}$ are mapped to.
This PDE in combination with periodic boundaries possesses an analytical solution of the form:
\begin{equation}
    \alpha(x,t) = \alpha(x-ut,0)
\end{equation}
That is, after traversing the quadratic domain with the given velocity $\underline{u} = \begin{bmatrix}1 & 1 \end{bmatrix}^T$, the solution field must exactly correspond to the initial condition.
Verification of numerical results is thus very straightforward.

We solve this problem using the Firedrake problem-solving environment along with the popular libraries PETSc and Scotch for efficient parallel computing \cite{FiredrakeUserManual,McRae2016,Homolya2016,Chaco95,PTSCOTCH,Rathgeber2016,dalcinParallelDistributedComputing2011,petsc-efficient,petsc-user-ref}.
As the Finite Volume Method can simply be understood as a Discontinuous Galerkin Method of polynomial degree zero, the implementation is virtually the same for both schemes.

Note that for the Finite Volume Method, the last term in Eq. \ref{eq:advection-weakform} becomes zero since the derivative of a constant vanishes. Thus, we omit this term from the assembly to save computations and to more accurately represent the arithmetic intensity posed by the original formulation of this scheme. 

For the sake of visualization, we project the solution onto a first-degree space with $H^1$ continuity.
The equations are solved by a three-stage implicit Runge Kutta method.
In both cases, careful attention has to be paid regarding the time step.
For hyperbolic problems of such time, the time step where stability is given is strictly bounded by the Courant Friedrichs Lewy number $\text{CFL} \leq \frac{1}{2k+1}$ for a scheme of degree $k$ \cite{cockburnRungeKuttaDiscontinuous2001}.
One can easily verify that for a Finite Volume scheme, this corresponds to the well-known condition that the CFL number must stay at or below unity.
Not only does that mean regarding arithmetic complexity that the FVM has a simplified assembly process, but also that the admissible time step is in general larger than for DG methods.
This discrepancy drastically increases with the polynomial order taken for the DGM.

\begin{figure}[htbp!]
    \centering
    \begin{subfigure}[t]{0.49\textwidth}
        \centering
        \includegraphics[width=\textwidth]{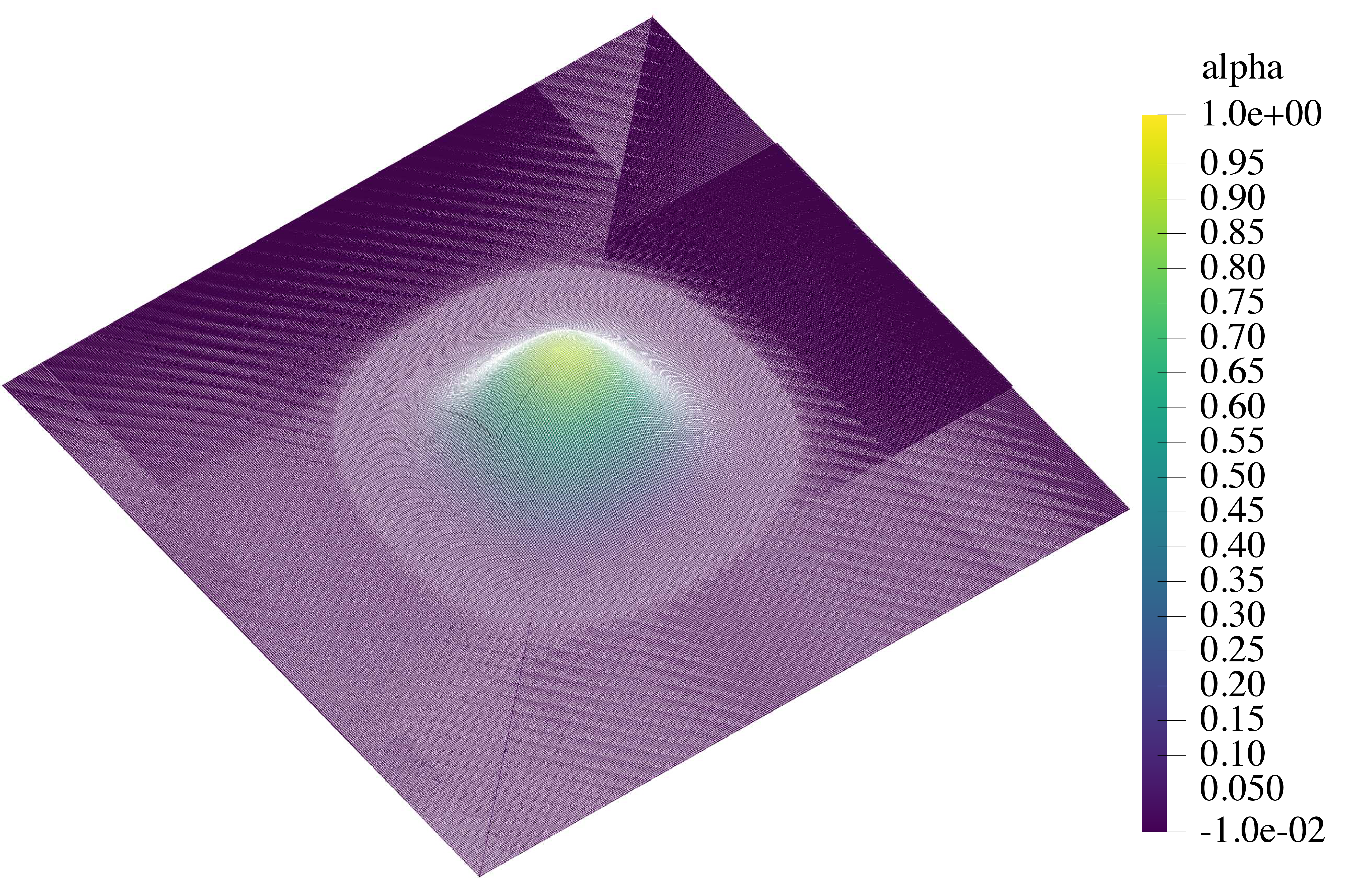}
        \caption{Finite Volume solution at end time on a 384x384 grid. The solution is heavily diffused to a parabolic profile.}
        \label{fig:advection-fvm}
    \end{subfigure}
    \hfill
    \begin{subfigure}[t]{0.49\textwidth}
        \centering
        \includegraphics[width=\textwidth]{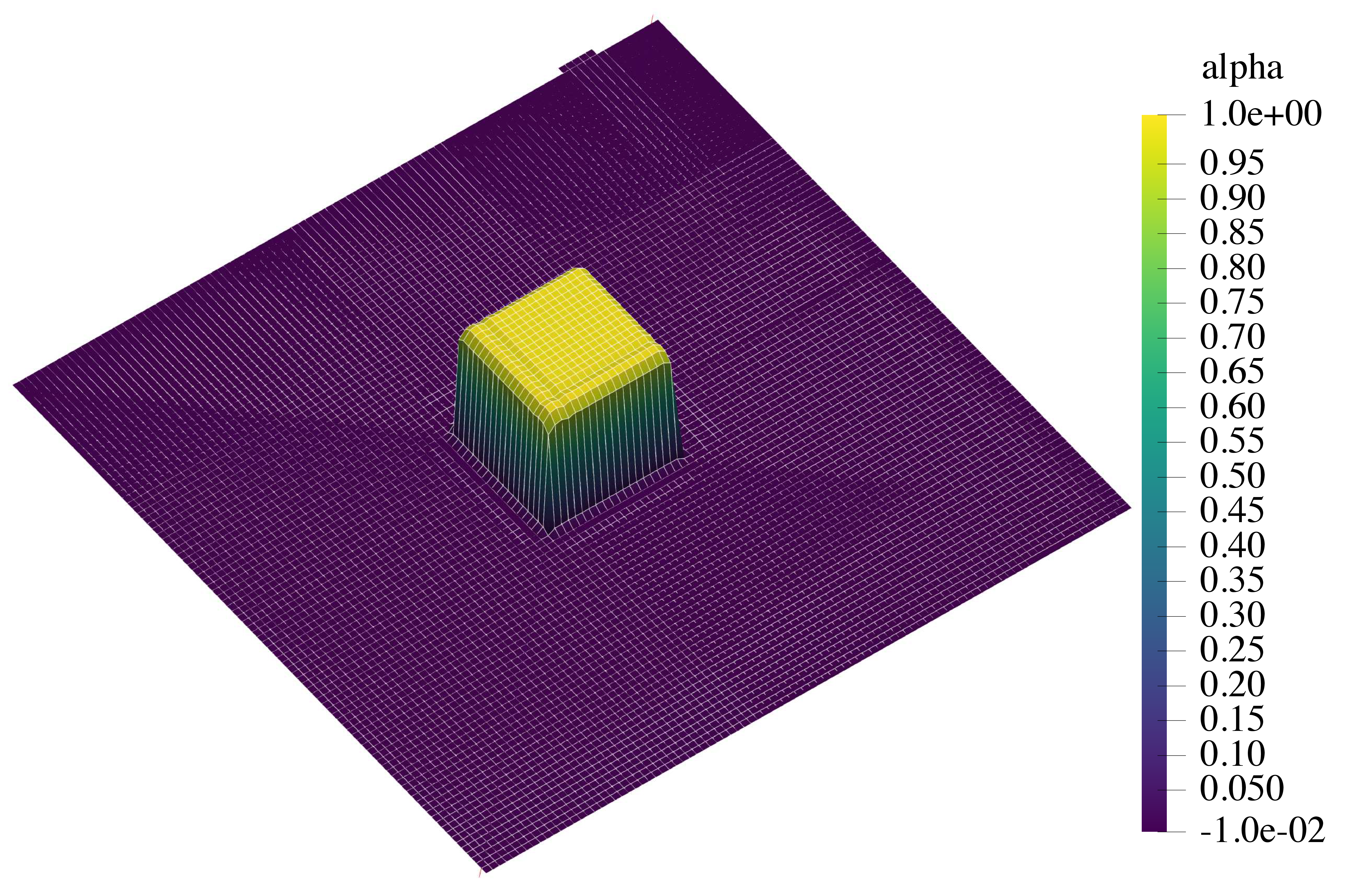}
        \caption{Discontinuous Galerkin solution at end time on a 96x96 grid using third order polynomials. One can observe the presence of spurious oscillations that remain stable in magnitude throughout the simulation.}
        \label{fig:advection-dgm}
    \end{subfigure}
   \caption{Solution of the two-dimensional advection equation using the Finite Volume and Discontinuous Galerkin method. Both models were run using the identical amount of global degrees of freedom set to 147.456. Element facets are drawn in white to illustrate the difference in grid size.}
   \label{fig:advection-solutions}
\end{figure}

The corresponding results of the simulation are shown in Figure \ref{fig:advection-solutions}.

The differences in accuracy are obvious.
The DG method can capture the rectangular profile throughout the simulation with relatively good accuracy, while the FV simulation is strongly diffused.

We note at this point that there are formulations of the FVM that capture shocks much more accurately whilst controlling oscillations.
Due to the maturity of this method, the field of constructing Weighted Essentially Non-Oscilaltory schemes (WENO) is quite advanced and also, in this case, will yield better approximations.
Such schemes are also applicable to considerably more complex systems of equations \cite{zimbrodEfficientSimulationComplex2022}.
However, this argument also applies to the DGM since it is as has been shown an extension of the FVM.
To make the comparison fair and to be able to rely on existing tools, we chose to compare both schemes with the same, relatively simple reconstruction technique.
For both cases, the use of WENO schemes would increase the computational load considerably, as these need to reconstruct polynomial approximations for each flux using relatively wide local stencils - similar to a high order Finite Difference scheme \cite{liuRobustReconstructionUnstructured2013}.
The choice of higher-order reconstruction techniques however should not affect the qualitative difference in approximation properties and performance.
For example, Zhou et al. report very similar findings for higher-order WENO schemes for the advection equation in two dimensions \cite{zhouNumericalComparisonWENO2001}.

\begin{figure}[htbp!]
    \centering
    \captionsetup[subfigure]{justification=centering}
    \begin{subfigure}[c]{0.49\textwidth}
        \centering
        \includegraphics[]{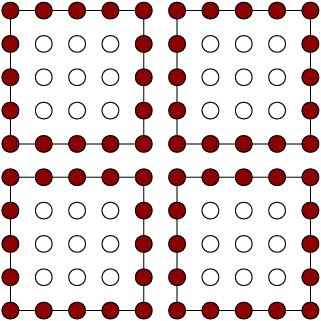}
        \caption{}
        \label{fig:dg-2d-shared_dofs}
    \end{subfigure}
    \hfill
    \begin{subfigure}[c]{0.49\textwidth}
        \centering
        \includegraphics[]{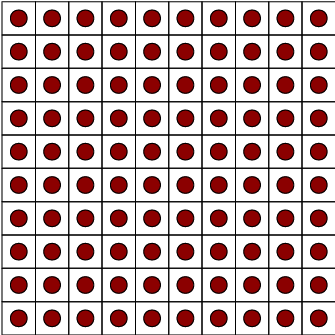}
        \caption{}
        \label{fig:fv-2d-shared_dofs}
    \end{subfigure}
   \caption{Comparison of shared DoFs (colored dots) between cells for a fourth-order DG method with an FV method with an equal amount of total DoFs (100).
   As the FVM is restricted to one DoF per cell, there are no interior DoFs, increasing the necessary MPI effort.
   The DGM on the other hand has an increasing amount of interior DoFs that for a collocation method do not contribute to the numerical flux.
   For the present order four scheme, the ratio of interior to total DoFs is 36\%.}
   \label{fig:fv-dg-shared_dofs}
\end{figure}

In addition to the previous comparison regarding accuracy, we also benchmark the capabilities for parallel computing using the three machines given in Table \ref{tab:advection-hardware}.
To properly scale up this problem, we aim to keep the number of Degrees of Freedom per core constant for the three environments, that is,  in this case around 18.400 DoFs/Core.
This is in the vicinity of the 25.000 DoFs/Core regime, where beyond that point considerable drop-offs in performance are to be expected \cite{badiaGenericFiniteElement2020}.
According to the theory, the DGM should perform with a noticeably higher efficiency due to less communication overhead between processes.
This is due to the reduced amount of cells for the same amount of DoFs and as such, the DGM has a much higher ratio of DoFs that lie inside the cell instead of at the boundary. 
As a consequence, there are fewer DoFs in relative terms that require the evaluation of a numerical flux and thus not as much overhead due to MPI efforts.
This relationship is visualized for reference in Figure \ref{fig:fv-dg-shared_dofs}.

The benchmarking results for all three hardware configurations are given in Figure \ref{fig:advection-runtimes}.
We report the total run times as well as the solution time for one singular time step.
This is since as explained the DG model has a considerably smaller admissible time step.
Thus, comparing overall run times is not suitable to empirically validate the above claims regarding parallel efficiency.
The overall computation time is nonetheless an important factor in terms of practicality since this is the main quantity one is interested in when performing a simulation.
\begin{figure}[htbp!]
\begin{minipage}[]{0.63\textwidth}
\centering
\input{Figs/Examples/Advection/advection.pgf}
\end{minipage}
\begin{minipage}[]{0.37\textwidth}
\centering
\begin{tblr}{
    cell{1}{3} = {c=2}{c},
    cell{2}{3} = {t},
    cell{2}{4} = {t},
    cell{3}{1} = {r},
    cell{4}{1} = {r},
    cell{5}{1} = {r},
    hline{1,6} = {-}{0.08em},
    hline{2} = {3-4}{},
    hline{3} = {-}{},
  }
   &  & \textbf{Wall time [s]} & \\
  \textbf{Cores} &  & \textbf{FVM} & \textbf{DGM}\\
  8 &  & 45.01 & 77.76\\
  18 &  & 114.77 & 136.92\\
  128 &  & 409.12 & 434.22
  \end{tblr}
\end{minipage}
\caption{Comparison of run times for the DGM and FVM advection benchmark case, run on the three configurations given in table \ref{tab:advection-hardware}.\\
\textit{Left:} Weak scaling of both models, measured as wall time per time step solve.
The amount of DoFs/Core was fixed at 18.400.
For each sample population, $n \approx 3000$.
Horizontal lines: Median, Boxes: IQR, Whiskers: Quartiles ± 1.5 IQR, \textit{****}: $p < 0.0001$.\\
\textit{Right:} Total Wall times for the solution over all time steps.
}
\label{fig:advection-runtimes}
\end{figure}
It becomes evident from the reported wall times that the FV model only has an advantage on the small desktop machine concerning solution time.
As predicted, scaling up the problem size and number of workers will yield an increasing advantage for the higher-order DGM.
The differences in computation time which were found to be empirically significant grow with increasing problem size, thus supporting the claim that DG methods are more favorable for highly parallel architectures.

A similar trend is visible from the table reporting overall runtimes.
The time step restrictions as discussed weaken the computational advantage gained in the solution of the semidiscrete system.
However, the difference in run time still decreases noticeably with increasing problem and hardware size.
For the largest machine with 128 cores, solution times are almost comparable, whereas in the case of the smallest machine, the FV model computed a solution about 42\% faster.

\subsection{Laser Powder Bed Fusion}

As a last example, we will investigate a real multiphysics problem to demonstrate the proposed method for more complex settings.
We outline the results from the problem analysis summarised in Figure \ref{fig:decision-process} and indicate how one may implement the resulting numerical scheme.

The problem of interest is the fluid flow and temperature field that is present during the laser-based heating of metallic powder.
Combined with a moving heat source that travels along a set path and layer-wise re-application of fresh powder, this technique results in the additive manufacturing process of Powder Bed Fusion of metallic materials using a laser beam (PBF-LB/M).

Having detailed knowledge about the fluid flow and temperature fields is crucial to obtain dense, that is, pore-free parts with a favorable microstructure.
The former is governed by the morphology of the solidified melt pool. The latter is heavily influenced by the spatial and temporal temperature gradients.

This specific problem has been tackled many times in the literature using various numerical methods.
An encompassing review of the technology as well as recent simulation approaches can be found in \cite{chowdhuryLaserPowderBed2022}

\subsubsection{Physics and Governing Equations}

Due to the presence of solid, liquid and gaseous phases, resolving the thermo fluid dynamics of the process involves a multitude of physics.
Most prominently, research has shown that the accurate representation of surface tension forces due to temperature gradient plays a key role in obtaining realistic results regarding the morphology of the solidified tracks, temperature gradients and presence of pores - all of which have shown to be important indicators for process quality \cite{debroyAdditiveManufacturingMetallic2018}.

The multiphase problem that arises from this application is discretized using the Volume of Fluid method.
The equations for this model have been introduced in section \ref{sec:advection-equation}.
In total, we have one solid metallic, one liquid metallic and two gaseous phases for vaporized metal and the shielding gas, yielding four phase fractions that need to be tracked.
Following the previous discussion in section \ref{sec:advection-equation}, one may omit to model one of these phases using a conservation law as one can calculate it using the compatibility condition in Eq. \ref{eq:vof-sum} as well.

With the phase fractions in place, one then obtains the mixed material constants at a given degree of freedom by a rule of mixture law.
For thermal conductivity $\kappa$, for instance, the phase averaged value at DoF $i$ is:
\begin{equation}
    \kappa_\text{VOF} = \sum\limits_{k=0}^{N_\text{Phases}} \kappa_k \alpha_k
\end{equation}
We continue by enumerating the governing equations of this problem.
For the melting of metallic materials, one finds that the flow field can be described by the incompressible Navier Stokes equations in the low Reynold's number regime, plus some source terms that account for the additional physics.
\begin{equation}\label{eq:momentum-balance}
    \partial_t (\rho \underline{u}) + \nabla \cdot \left[\left(p - p_\text{recoil}\right) \underline{\underline{I}}\right] - \rho \underline{k} - \eta \Delta \underline{u} + \nabla \cdot \underline{\underline{T}}_\text{capillary}  = 0
\end{equation}
With the additional condition for incompressibility:
\begin{equation}\label{eq:continuity}
    \nabla \cdot \underline{u} = 0
\end{equation}
Here, $\underline{k}$ is the vector of gravity and other volumetric forces.

One important driving force for melt flow is the temperature-dependent surface tension force, also called the Marangoni force.
We can summarise all capillary forces normal and tangential to the interface using the so-called capillary stress tensor:
\begin{equation}\label{eq:capillary-stress}
    \underline{\underline{T}}_\text{capillary} = - \nabla \cdot \left[-\sigma (\underline{\underline{I}} - \underline{n} \otimes \underline{n}) \right]
\end{equation}
Where $\sigma$ denotes the temperature-dependent coefficient of surface tension, $I$ is the identity tensor and $n$ is the unit normal vector to the capillary surface.
The divergence of this quantity represents the capillary force acting on interfaces \cite{brackbillContinuumMethodModeling1992}.

This in turn needs to be computed using the phase fractions and the surface gradient operator $\nabla_s$:
\begin{align}
    \underline{n} = \frac{\nabla_s \alpha}{\lvert \nabla_s \alpha \rvert} \label{eq:unit-normal} \\
    \nabla_s(\alpha_1,\alpha_2) = \alpha_1 \nabla \alpha_2 - \alpha_2 \nabla \alpha_1 \label{eq:interface-gradient}
\end{align}
For the gaseous phase of vaporized material, we additionally consider the recoil pressure
\begin{equation}\label{eq:recoil-pressure}
    p_\text{recoil}(T) = 0.53 p_0 \exp{\frac{L_v}{R} \left( \frac{1}{T_v}-\frac{1}{T} \right)}
\end{equation}
The only conserved quantity left is the total energy of the system.
The corresponding balance equation that governs the evolution of temperature reads:
\begin{equation}\label{eq:heat-balance}
    \partial_t (\rho c_p T) + \underline{u} \cdot \nabla (\rho c_p T) + \sum_i \rho L \partial_t \alpha_i - \nabla \cdot (\kappa \nabla T)  + Q_\text{Laser} - Q_\text{Vap} - Q_\text{Radiation} = 0
\end{equation}
With the following (nonlinear) source terms for the laser heat input, vaporization loss and radiation loss, respectively:
\begin{align}
    Q_\text{Laser}(x,y) = \frac{2P}{\pi r_0^2} \exp \left( - \frac{2(x^2+y^2)}{r_0^2}\right) \label{eq:laser-heat}\\
    Q_\text{Vap}(T) = 0.82 \frac{p_\text{recoil}}{\sqrt{2\pi MRT}} \label{eq:vapor-heat} \\
    Q_\text{Radiation}(T) = \sigma \epsilon (T^4 - T^\text{amb}) \label{eq:radiation-heat}
\end{align}
The laser that provides $Q_\text{Laser}$ travels along a straight line throughout the domain with constant velocity and laser power $P$ with beam radius $r_0$.
Thus, this source term represents a transient heat source where its position can be easily interpolated given the current time.
Further details on the physical models can, for example,  be found in \cite{zimbrodEfficientSimulationComplex2022}.

\subsubsection{Computational Domain}

We solve the given problem on a grid that is overall three-dimensional and hexahedral in shape.
The bounding box has dimensions \SI{2.0}{\milli\metre} in length, \SI{0.3}{\milli\metre} in width and \SI{0.5}{\milli\metre} in height.

In this problem setting, steep temperature gradients are to be expected due to the concentrated heat input.
As a result, it makes sense to employ a finer grid at locations where the melt pool dynamics take place.
One possible way to account for that is to use a discretization technique that can incorporate adaptive grids.

However, we would instead like to make use of our previous knowledge and only use a finer grid around the vicinity of the laser path - where the actual melting takes place and a finer grid is needed to appropriately resolve the flow field.
As such, there is no need to use algorithms that flag candidate cells and employ re-meshing at every time step.
The rest of the domain can be meshed using coarser cells, as they are primarily present in the simulation to properly account for heat conduction to surrounding solidified material.

\begin{figure}[htbp!]
    \centering
    \includegraphics[width=0.6\linewidth]{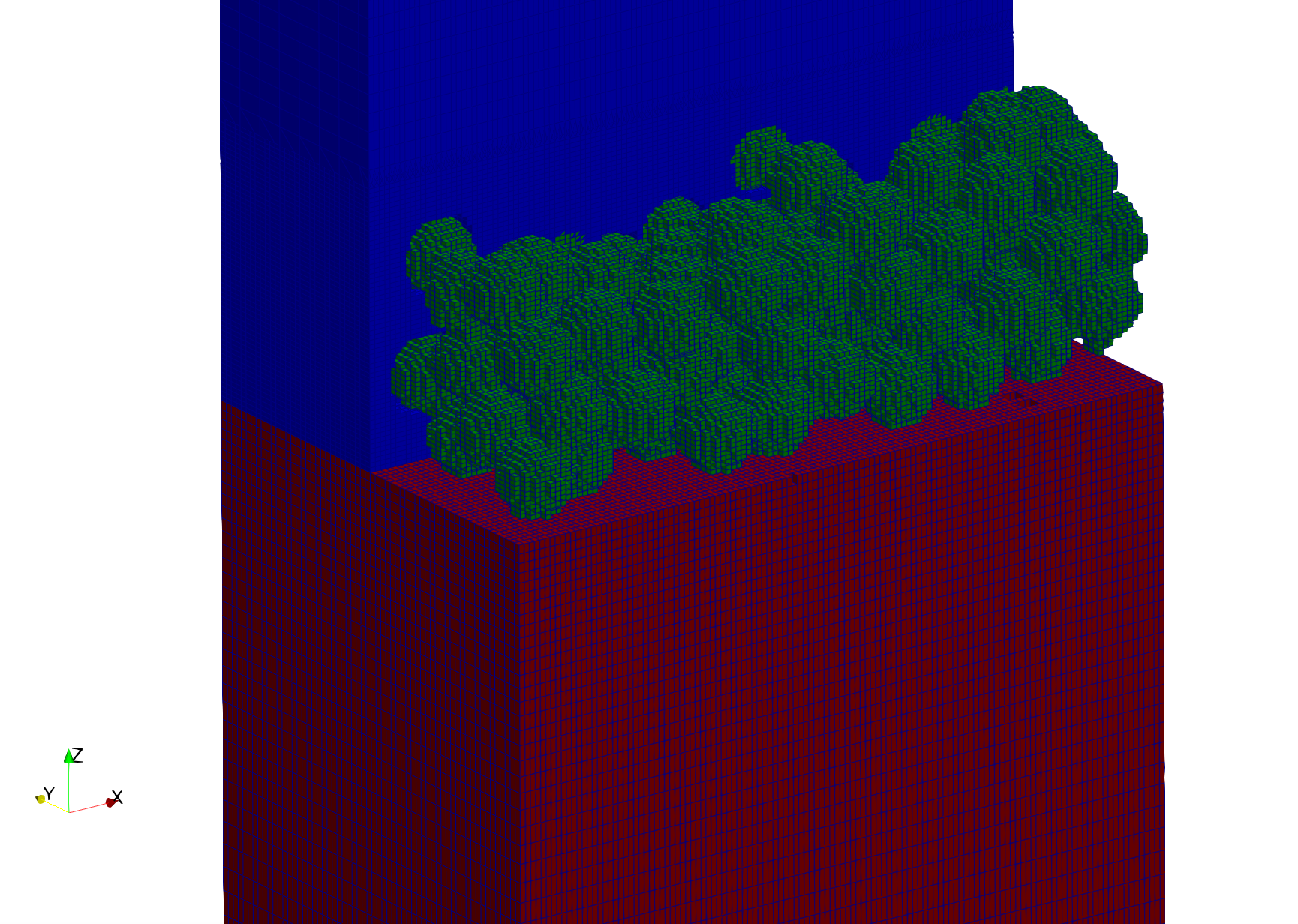}
    \caption{Slice of the computational domain for the Laser Powder Bed Fusion test case. The grid consists of hexahedral elements with varying sizing. Colors denote different phases in the initial condition. This figure is adapted from \cite{zimbrodEfficientSimulationComplex2022}.}
    \label{fig:lpbf-mesh}
\end{figure}

A slice of the resulting discretization is depicted in Figure \ref{fig:lpbf-mesh}.
The actual mesh is box-shaped, and the visible regions are colored to illustrate the different phases present as well as the varying cell sizes.
One can observe that the hexahedral cells become coarser with increasing depth of the powder bed.

\subsubsection{Computing Resources}

We start the selection process for suitable numerical schemes as done in the previous sections by investigating the computational hardware, as stated within process P1 in Figure \ref{fig:decision-process}.
For this problem, we use the same hardware given in Table \ref{tab:advection-hardware}.
That is, we aim to solve this problem on the 18-core CPU machine.
As given in section \ref{sec:hardware}, such a system of desktop scale does not per se necessitate the use of massively parallel hardware, as the amount of processes falls below the recommended amount of 50 to 100 workers.

\subsubsection{Problem Scale}

Using the system of PDEs as input, we can proceed with analyzing the problem length scales according to process P2 in Figure \ref{fig:decision-process}.
Given the coefficients and material properties of the PDE system, we note that all physics are expected to operate on a length scale of a few micrometers.
We can thus conclude that this particular problem does not exhibit multiscale properties that would necessitate adaptivity.
Therefore, in combination with the given hardware, we can negate the first decision point D1 in Figure \ref{fig:decision-process} and proceed with the field-wise analysis of the PDE system.

\subsubsection{Classification}

We begin with the field-wise classification of the PDE system, that is,  decisions D2 to D5 along with the corresponding processes as follows:
We consider all governing equations of the system (Eqs. \ref{eq:capillary-stress}-\ref{eq:radiation-heat}) individually for each field quantity. That is, we extract the individual equations that one needs to solve.
Afterward, we follow the classification methods described by P3 and P4, whereas the latter is only relevant for hyperbolic systems as illustrated in Figure \ref{fig:decision-process}.

We first examine the phase fraction fields.
As all variables are governed by the same set of equations and thus the resulting classification applies to every field, we abbreviate the discussion by only classifying one phase fraction $\alpha$.
The relevant governing equation then is exactly the passive advection equation that has been introduced in section \ref{sec:advection-equation}:

\begin{equation}
    \partial_t \alpha + \underline{u} \cdot \nabla \alpha = 0
\end{equation}

It then immediately follows from the previous discussion that the phases are governed by first-order, linear, hyperbolic systems.
We consequently arrive at the Discontinuous Galerkin Method for these variables.

Next, we will deal with the velocity vector $\underline{u}$.
Extracting all components from Eq. \ref{eq:momentum-balance} that contain differentials of $\underline{u}$, we have:
\begin{equation}
    \rho \partial_t \underline{u} - \eta \Delta \underline{u} = 0
\end{equation}
This system of equations, for each component, possesses the exact structure of the heat diffusion equation, which is known to be a model parabolic equation.
Following the more rigorous approach given in section \ref{sec:pde-classification}, we would obtain for the coefficient matrix $a_{ij}$:
\begin{equation}
    a_{ij}^{\underline{u}} = \begin{bmatrix}
        0 & 0 & 0 \\
        0 & -\eta & 0 \\
        0 & 0 & -\eta
    \end{bmatrix}
\end{equation}
The first zero diagonal element is due to the time component only appearing in a first, but not in any second derivative.
We consequently arrive at the same conclusion, being that the governing equation for velocity in this case is parabolic.
Thus, we follow the left branch after D2 and continue by classifying the regularity of the domain.
The question of whether the grid is strictly regular and cartesian has already been answered in the previous section, and therefore we arrive at the Continous Galerkin Method for the velocity field.

For the pressure $p$, we find that this quantity only appears within the momentum balance in Eq. \ref{eq:momentum-balance} within the expression $\nabla \cdot p\underline{\underline{I}}$, which is equal to the gradient $\nabla p$.
That is, we have $\partial_x p$ as a source term in the components in this equation and thus the evolution of $p$ is not per se given by a separate governing equation.
We therefore fall back to the qualitative classification that was described in section \ref{sec:pde-classification} and choose an appropriate class based on the continuity requirements of the pressure field.
The present problem consists of multiple phases, therefore we may expect steep jumps in pressure at interfaces.
We consequently assume the pressure field to be discontinuous at some points during the simulation and thus follow the process in Figure \ref{fig:decision-process} along the path of hyperbolic equations.
As there are no additional source terms or nonlinearities in $p$, we classify this constraint as linear in D4 and therefore also arrive at the Discontinuous Galerkin Method.

For the remaining Temperature field $T$, we again gather the terms containing differential operators on $T$ from the governing equation, Eq. \ref{eq:heat-balance}:
\begin{equation}
    \rho c_p \partial_t T + \rho c_p \underline{u} \cdot \nabla T - \kappa \Delta T = \text{RHS}(T)
\end{equation}
The remaining terms in Eq. \ref{eq:heat-balance} that have been gathered in $\text{RHS}(T)$ are either constant, or depend on $T$ itself, but not in its derivatives. Therefore, we may omit them for the sake of classification.
As such, we have a PDE of second order and proceed analogously to the classification of the velocity $\underline{u}$.
The coefficient matrix $a_{ij}$ is then:
\begin{equation}
    a_{ij}^{T} = \begin{bmatrix}
        0 & 0 & 0 \\
        0 & -\kappa & 0 \\
        0 & 0 & -\kappa
    \end{bmatrix}
\end{equation}
Similarly to the previous discussion, this resembles a parabolic system and we thus recommend the discretization of the temperature field using the Continuous Galerkin Method, as the question of domain regularity (D3) has already been addressed and is the same for all fields.
If the coefficient of thermal diffusion $\kappa$ would be zero in this case, the governing equation would then collapse into a hyperbolic PDE and another classification would apply.
We can however state that the amount of heat conduction present in the model cannot be neglected, and thus we have that $\kappa > 0$ everywhere in the domain.

\subsection{Resulting discretization}

With all independent variables of the model addressed in the previous section, we briefly summarise the chosen combination of numerical schemes in Table \ref{tab:lpbf-summary}.
\begin{table}
    \centering
    \caption{Resulting schemes of the proposed decision method for the Laser Powder Bed Fusion example problem. Columns D1 to D4 refer to the corresponding decision points in Figure \ref{fig:decision-process}.}
    \label{tab:lpbf-summary}
    \begin{tabular}{llllll} 
    \toprule
    \textbf{Variable} & \textbf{D1} & \textbf{D2} & \textbf{D3} & \textbf{D4} & \textbf{Result} \\ 
    \hline
    $\alpha_\text{solid}$ & no & hyperbolic & n.a. & linear & DGM \\
    $\alpha_\text{liquid}$ & no & hyperbolic & n.a. & linear & DGM \\
    $\alpha_\text{gas}$ & no & hyperbolic & n.a. & linear & DGM \\
    $\underline{u}$ & no & parabolic & irregular & n.a. & CGM \\
    $p$ & no & hyperbolic & n.a. & linear & DGM \\
    $T$ & no & parabolic & irregular & n.a. & CGM \\
    \bottomrule
    \end{tabular}
    \end{table}
One may proceed to implement this mixed discretization using any modern Finite Element code.
In the first two examples, we have shown that both creating a custom implementation using modern programming languages and implementing a model using popular Finite Element libraries are viable options to create a corresponding simulation.

In this case, we obtain a coupled, nonlinear problem.
This means that instead of solving a discrete linear system of equations, one must solve a global, nonlinear root finding problem of the form $F(u;v) = 0$ in every time step by appropriate means, such as Picard iteration or Newton's method.
Such functionality however is readily implemented in Finite Element libraries such as Firedrake \cite{FiredrakeUserManual}, which has been used in the previous section.

As a final remark, we note that similar models for resolving the melt pool dynamics in Laser Powder Bed fusion have been developed using mixed formulations of the Finite Element Method.
For example, Meier et al. developed a custom, high-performance Finite Element simulation code called \texttt{MeltPoolDG} based on the Discontinuous Galerkin Method as well as a mixed FEM to resolve the melt pool on the powder bed scale \cite{meierPhysicsBasedModelingPredictive2021,kronbichlerFastMassivelyParallel2018}.
Other recent developments include a Finite Element model on multiple grids to capture the governing physics \cite{caboussatNumericalSimulationTemperaturedriven2023}, and a space-time FEM code that is based on a Petrov Galerkin approximation \cite{koppSpacetimeHpfiniteElements2022}.
Furthermore, there also exist implementations of such melt pool models using commercially available software.
For instance, the popular COMSOL software package which implements many variants of the FEM can resolve the relevant thermo fluid dynamics during laser melting \cite{liuSimulationPowderPacking2020,mayiTransientDynamicsStability2021}.

\section{Conclusion and Future Work}

We presented a general framework for determining combinations of interoperable, grid-based approximation schemes that are individually suited for the solution of multiphysics problems.
The key findings of this work can be summarised as follows:
\begin{itemize}
    \item The Discontinuous Bubnov Galerkin Method can serve as a general mathematical baseline for the Finite Volume Method, Finite Difference Method and Continuous Galerkin Method.
    These form the majority of schemes used in practice, which is why we restrict the scope of this work to this specific subset.
    We have shown that one can recover all these schemes by introducing certain restrictions to the DG Method, which can either be on the algebraic or algorithmic level.
    \item By using those simplified schemes where appropriate, one is to gain stability and performance in numerical simulations that scale with problem size and hardware.
    In particular, one may even avoid the assembly of a global linear system by leveraging domain restrictions in the case of the Finite Difference Method.
    If a PDE exhibits strong nonlinearity, choosing a purely reconstruction-based approach via the Finite Volume Method is then beneficial and assembly of the weak form can be avoided as well.
    \item This common framework based on the DG method enables interoperability of the mentioned schemes.
    As such, one can combine them when solving multiphysics problems and thus assemble efficient mixed discretizations.
    Therefore, manual coupling of different numerical solvers is not needed in theory.
    Instead, one may implement them in a monolithic way and thus avoid costly memory transfer operations.
    \item One can systematically derive an efficient combination of numerical schemes, given some inputs about the problem, that is, Hardware requirements, the mathematical formulation of the problem and the domain geometry. The method to identify these schemes is based on the field and delivers one-to-one recommendations on which method to use.
    \item Applying the developed method to the two model problems in this work, the Allen Cahn equation in 1D and 2D and the advection equation in 2D, yields methods that notably outperform their respective alternatives. In the former case, this difference in computation time exceeds one order of magnitude.
    In the latter case, we have shown that as predicted by the theory, the computational advantages scale with problem size and the necessary degrees of freedom for some fixed accuracy differ considerably.
    We have also shown that for a real multiphysics problem, even for coupled systems of PDEs, one can derive an efficient combination of schemes in a reproducible and for application experts accessible way.
    This highlights the practicality and usefulness of the established framework for end users.
    Especially researchers that are familiar with physical modeling, but not to the same degree with state-of-the-art high-performance numerical methods are expected to profit from this approach.
\end{itemize}
Some future work might be attributed to further breaking down the choices of numerical schemes based on a more granular problem classification.
The classes of PDE we introduced in this work are rather coarse, albeit the most widely used.

One may also include mixed Finite Element discretizations into the method.
Within the scope of this work, we only included scalar- and vector-valued Lagrange elements.
Including mixed Finite Elements might however make the selection process ambiguous and thus special care must be attributed to keep the simplicity of the method.

%%%%%%%%%%%%%%%%%%%%%%%%%%%%%%%%%%%%%%%%%%
\vspace{6pt} 

%%%%%%%%%%%%%%%%%%%%%%%%%%%%%%%%%%%%%%%%%%
%% optional
%\supplementary{The following supporting information can be downloaded at:  \linksupplementary{s1}, Figure S1: title; Table S1: title; Video S1: title.}


%\begin{tabularx}{\textwidth}{lll}
%\toprule
%\textbf{Name} & \textbf{Type} & \textbf{Description} \\
%\midrule
%S1 & Python script (.py) & Script of python source code used in XX \\
%S2 & Text (.txt) & Script of modelling code used to make Figure X \\
%S3 & Text (.txt) & Raw data from experiment X \\
%S4 & Video (.mp4) & Video demonstrating the hardware in use \\
%\bottomrule
%\end{tabularx}
%}

%%%%%%%%%%%%%%%%%%%%%%%%%%%%%%%%%%%%%%%%%%
\authorcontributions{Conceptualization, P.Z. and M.F.; method, P.Z.; software, P.Z.; validation, P.Z. and M.F.; formal analysis, P.Z.; investigation, P.Z. and M.F.; resources, J.S.; writing---original draft preparation, P.Z.; writing---review and editing, M.F. and J.S.; visualization, P.Z. and M.F.; supervision, J.S. All authors have read and agreed to the published version of the manuscript.}

\dataavailability{All data supporting this publication, including the manuscript source, are available under \href{https://github.com/pzimbrod/multiphysics-pde-methods}{https://github.com/pzimbrod/multiphysics-pde-methods}.} 

\conflictsofinterest{The authors declare no conflict of interest.} 

%%%%%%%%%%%%%%%%%%%%%%%%%%%%%%%%%%%%%%%%%%
%% Optional

%% Only for journal Encyclopedia
%\entrylink{The Link to this entry published on the encyclopedia platform.}

\abbreviations{Abbreviations}{
The following abbreviations are used in this manuscript:\\

\noindent 
\begin{tabular}{@{}ll}
CFL & Courant-Friedrichs-Lewy \\
CG(M) & Continuous Galerkin (Method) \\
DG(M) & Discontinuous Galerkin (Method) \\
DoF & Degree of Freedom \\
FD(M) & Finite Difference (Method) \\
FV(M) & Finite Volume (Method) \\
MPI & Message Passing Interface \\
PDE & Partial Differential Equation \\
VoF & Volume of Fluid
\end{tabular}
}

%%%%%%%%%%%%%%%%%%%%%%%%%%%%%%%%%%%%%%%%%%
%% Optional
\appendixtitles{yes} % Leave argument "no" if all appendix headings stay EMPTY (then no dot is printed after "Appendix A"). If the appendix sections contain a heading then change the argument to "yes".
\appendixstart
\appendix
\section[\appendixname~\thesection]{Implementation Details}
\subsection[\appendixname~\thesubsection]{Allen Cahn Equation}

We solve both systems of Ordinary Differential Equations (ODE) using the \texttt{DifferentialEquations.jl} library in Julia.
Details about model parametrization are given in Table \ref{tab:ac-parameters}.

\begin{table}[htbp!]
    \centering
    \caption{Parameters used for both 1D and 2D models of the Allen Cahn equation}
    \label{tab:ac-parameters}
    \begin{tabular}{llll} 
    \toprule
    \textbf{Variable} & \textbf{1D Model} & \textbf{2D Model} & \textbf{Description} \\ 
    \hline
    $\Gamma$ & 1.0 & 50.0 & Interface energy \\
    $\xi$ & 1.5 & 4.0 & Interface width \\
    $M$ & 1.0 & 1.0 & Kinetic mobility \\
    $\mu_0$ & 0.1 & 0.0 & Bulk energy density differential \\
    $x_0$ & 20.0 & n.a. & Initial interface position \\
    $h$ & 1.0 & 1.0 & Grid size \\
    \bottomrule
    \end{tabular}
    \end{table}

We benchmark both models using the \texttt{BenchmarkTools.jl} library.
The physical size of both models was measured using the Julia Base function \texttt{summarysize}.
For timing purposes, each run was repeated 100 times to create reliable empirical data.
Both models were run on the same machine in consecutive runs.

\subsection[\appendixname~\thesubsection]{Advection Equation}

\subsubsection[\appendixname~\thesubsubsection]{Derivation of the Weak Form}

Eq. \ref{eq:advection-weakform} can be obtained by the strong form of the PDE as described in the article, by applying partial integration after multiplying with a test function and integrating over the domain $\Omega$.
To arrive at a formulation that can be implemented, there are a few additional steps which are described below.

In most DG codes, the weak form given by Eq. \ref{eq:advection-weakform}, more specifically the surface integral cannot be entered directly but needs to be reformulated using jump and average operators.
This is also done in the advection equation example described in section \ref{sec:advection-equation}.
Such a modified weak form can be derived as follows.

The most straightforward consists of summarizing the convective term of the PDE as the divergence of some flux $\nabla \cdot \underline{f}(\alpha)$.
In this case, the flux is simply $\underline{f}(\alpha) = \underline{u} \phi$.

Then, the surface integral of Eq. \ref{eq:advection-weakform} can be formulated as:
\begin{equation}
    \int_\Gamma v \alpha (\underline{u}\cdot\underline{n})\,dS = \int_\Gamma v \underline{f}(\alpha)\cdot\underline{n}\,dS
\end{equation}
We now write down the terms of the surface integral of Eq. \ref{eq:advection-weakform} that contribute to an arbitrary internal facet, that is,  we always encounter two integrals from the cells that the facet is owned by.
We distinguish between the two cells by introducing the notation "$+$" and "$-$":
\begin{equation}\label{eq:advection-facet}
    \int_\Gamma v_h^+ \alpha_h^+ (\underline{u}^+\cdot\underline{n}^+) \,dS + \int_\Gamma v_h^- \alpha_h^- (\underline{u}^-\cdot\underline{n}^-) \,dS
\end{equation}
Where we encounter the yet-to-be-determined values for $\alpha^+$ and $\alpha^-$.
By introducing the common, yet unknown numerical flux $\underline{f}(\alpha)^*$ and utilizing the fact that the facet normal vectors $\underline{n}$ point opposite to each other, Eq. \ref{eq:advection-facet} becomes:
\begin{equation}
     \int_\Gamma \left(\underline{f}(\alpha)^* \cdot \underline{n}\right)(v^+ - v^-) \,dS = \int_\Gamma \left(\underline{f}(\alpha)^* \cdot \underline{n}\right) [\![v]\!] \,dS
\end{equation}
The only issue that is left to deal with is now the choice of an appropriate expression for the numerical flux $\underline{f}(\alpha)^* \cdot \underline{n}$.

Both DG and FV implementations of the advection equation use the Lax-Friedrichs flux function:
\begin{equation}
    f(\alpha)^* = \{ \alpha (\underline{u} \cdot \underline{n}^+) \} + \frac{1}{2} C [\![ \alpha ]\!]  
\end{equation}
Where $\{\cdot\}$ and $[\![\cdot]\!]$ denote the average and jump operators as defined above.
In this case, $C$ is the maximum velocity in the domain and takes the role of the maximum signal velocity, which must be computed separately for more complex problems.
It is defined as the jacobian of the PDE flux $\partial f(\alpha) / \partial \alpha$.
Since we have that $f(\alpha) = \underline{u} \alpha$, the jacobian simply becomes the prescribed advection velocity $\underline{u}$.

\subsubsection[\appendixname~\thesubsubsection]{Simulation Parameters}

We solve the resulting semi-discrete ODE systems using a three-stage implicit Runge Kutta formulation using a matrix-free solver.

We time the entire simulation runs from the command line and execute each run six times.
Wall times per time step were recorded during the run using the PETSc distributed interface and logged to a text file.

%%%%%%%%%%%%%%%%%%%%%%%%%%%%%%%%%%%%%%%%%%
\begin{adjustwidth}{-\extralength}{0cm}

\reftitle{References}


\bibliography{bibliography}

% If authors have biography, please use the format below
%\section*{Short Biography of Authors}
%\bio
%{\raisebox{-0.35cm}{\includegraphics[width=3.5cm,height=5.3cm,clip,keepaspectratio]{Definitions/author1.pdf}}}
%{\textbf{Firstname Lastname} Biography of first author}
%
%\bio
%{\raisebox{-0.35cm}{\includegraphics[width=3.5cm,height=5.3cm,clip,keepaspectratio]{Definitions/author2.jpg}}}
%{\textbf{Firstname Lastname} Biography of second author}

% For the MDPI journals use author-date citation, please follow the formatting guidelines on http://www.mdpi.com/authors/references
% To cite two works by the same author: \citeauthor{ref-journal-1a} (\citeyear{ref-journal-1a}, \citeyear{ref-journal-1b}). This produces: Whittaker (1967, 1975)
% To cite two works by the same author with specific pages: \citeauthor{ref-journal-3a} (\citeyear{ref-journal-3a}, p. 328; \citeyear{ref-journal-3b}, p.475). This produces: Wong (1999, p. 328; 2000, p. 475)

%%%%%%%%%%%%%%%%%%%%%%%%%%%%%%%%%%%%%%%%%%
%% for journal Sci
%\reviewreports{\\
%Reviewer 1 comments and authors’ response\\
%Reviewer 2 comments and authors’ response\\
%Reviewer 3 comments and authors’ response
%}
%%%%%%%%%%%%%%%%%%%%%%%%%%%%%%%%%%%%%%%%%%
\PublishersNote{}
\end{adjustwidth}
\end{document}